\documentstyle[11pt,a4paper,leqno,amscd,amssymb,pstricks,latexsym,amsbsy,xypic,mathrsfs,verbatim,tikz]{amsart}

\newtheorem{thm}{Theorem}[section]
\newtheorem{lem}[thm]{Lemma}

\newtheorem{prop}[thm]{Proposition}

\theoremstyle{definition}
\newtheorem{example}[thm]{Example}

\newtheorem{defn}[thm]{Definition}

\newtheorem{conj}[thm]{Conjecture}
\newtheorem{rem}[thm]{Remark}

\numberwithin{equation}{thm}

\input xy
\xyoption{all}

\begin{document}

\title[On tropical friezes associated with Dynkin diagrams]{On tropical friezes associated \\ with Dynkin diagrams}

\author{Lingyan GUO}
\address{Universit\'e Paris Diderot - Paris~7, UFR de Math\'ematiques,
Institut de Math\'ematiques de Jussieu, UMR 7586 du CNRS, Case 7012,
B\^atiment Chevaleret, 75205 Paris Cedex 13, France}
\email{guolingyan@@math.jussieu.fr}

\date{\today}

\begin{abstract}
Tropical friezes are the tropical analogues of Coxeter-Conway frieze
patterns. In this note, we study them using triangulated categories.
A tropical frieze on a $2$-Calabi-Yau triangulated category
$\mathcal{C}$ is a function satisfying a certain addition formula.
We show that when $\mathcal{C}$ is the cluster category of a Dynkin
quiver, the tropical friezes on ${\mathcal{C}}$ are in bijection
with the $n$-tuples in ${\mathbb{Z}}^n$, any tropical frieze $f$ on
$\mathcal{C}$ is of a special form, and there exists a
cluster-tilting object such that $f$ simultaneously takes
non-negative values or non-positive values on all its indecomposable
direct summands. Using similar techniques, we give a proof of a
conjecture of Ringel for cluster-additive functions on stable
translation quivers.

\end{abstract}

\maketitle

\section{introduction}

Cluster algebras introduced by S. Fomin and A. Zelevinsky \cite{FZI}, are subrings of the field ${\mathbb{Q}}(x_1, \ldots,x_n)$
of rational functions in $n$ indeterminates endowed with a distinguished set of generators called cluster variables, which are constructed
recursively via an operation called mutation. A cluster algebra is said to be of finite type if the number of cluster variables
is finite. The classification of finite type cluster algebras was achieved \cite{FZII} in terms of Dynkin diagrams.

Motivated by close relations between tilting theory of finite-dimensional hereditary algebras and the combinatorics of mutation in cluster
algebras, the cluster category ${\mathcal{C}}_Q$ of a finite acyclic quiver $Q$ was introduced in \cite{CCS} for type $A_n$ and in
\cite{BMRRT06} for the general case. The cluster category provides a natural model for the combinatorics of its corresponding cluster
algebra. It is triangulated \cite{Ke05}, Krull-Schmidt and 2-Calabi-Yau \cite{BMRRT06} in the sense that there are bifunctorial
isomorphisms $${\rm Ext}^1(X,Y) \simeq D {\rm Ext}^1(Y,X), \quad X, Y \in {\mathcal{C}}_Q.$$
There are also many other 2-Calabi-Yau triangulated categories, for example, the stable module categories of preprojective algebras
of Dynkin type studied by Geiss-Leclerc-Schr${\ddot{\mathrm{o}}}$er in their series of papers, the generalized cluster categories
of Jacobi-finite quivers with potential \cite{DWZ} and of finite-dimensional algebras of global dimension $\leq 2$, which were investigated
in \cite{Am08} by C. Amiot.

Starting from a 2-Calabi-Yau Hom-finite triangulated category $\mathcal{C}$ with a cluster-tilting object $T$, Palu \cite{Pal08} introduced
the notion of a cluster character $\chi$ from $\mathcal{C}$ to a commutative ring which satisfies the multiplication formula
$$\chi(L) \chi(M) = \chi(E) + \chi(E')$$ for all objects $L$ and $M$ such that ${\rm Ext}^1_{\mathcal{C}} (L,M)$ is one-dimensional,
where $E$ and $E'$ are the middle terms of the non-split triangles
%$$L \rightarrow E \rightarrow M \rightarrow \Sigma L \quad \mbox{and} \quad
%M\rightarrow E' \rightarrow L \rightarrow \Sigma M$$
with end terms $L$ and $M$.
%He showed that the map $L \mapsto X^T_L$ (a fraction
%using a formula proposed by Caldero-Keller \cite{CK08}) is a cluster character and `categorifies' the corresponding cluster algebra in the acyclic case.
He explicitly constructed cluster characters from cluster-tilting objects.

In this article, we introduce tropical friezes $f$ on $\mathcal{C}$
mainly by replacing the above multiplication formula with an
addition formula $$f(L) + f(M) = \max \{ f(E), f(E')\}.$$ Our
inspiration comes from the definition of cluster-additive functions
\cite{Rin11} on stable translation quivers and from the tropicalized
version of Coxeter-Conway's frieze patterns. To the best of our
knowledge, such tropical frieze patterns first appeared implicitly
in Fock-Goncharov's preprint \cite{FG1} and explicitly in section 4
of J. Propp's preprint \cite{JP}.

The paper is organized as follows.

In Section 2, after recalling some facts on frieze patterns and stating the assumptions on the categories $\mathcal{C}$ we consider (namely, 2-Calabi-Yau categories
with cluster-tilting object), we introduce the notion of tropical friezes. Then we study their first properties and some links
to cluster characters, using which we give an example and a counter-example of tropical friezes.

In Section 3, taking advantage of the indices \cite{KR07} of objects of ${\mathcal{C}}$, for each cluster-tilting object $T$ and
each element $m$ in the Grothendieck group $K_0({\rm mod} {\rm End}_{\mathcal{C}}(T))$, we define a function $f_{T,m}$ on $\mathcal{C}$.
A criterion for $f_{T,m}$ to be a tropical frieze is given in Theorem \ref{2}, which is also a necessary condition when $\mathcal{C}$
is the cluster category ${\mathcal{C}}_Q$ of a Dynkin quiver $Q$. We also show that the tropical friezes on ${\mathcal{C}}_Q$
with $Q$ Dynkin are in bijection with the $n$-tuples in ${\mathbb{Z}}^n$ by composing Palu's cluster character with a morphism of semifields. Then we investigate
the cluster-hammock functions introduced by Ringel \cite{Rin11}, which always give rise to tropical friezes while their sums do not.

Section 4 just consists of simple illustrations for the cases $A_1$ and $A_2$, in order to give the reader an intuitive impression.

In Section 5, for a cluster-tilting object $T$ and a tropical frieze $f$ on ${\mathcal{C}}$, we define an element $g(T)$ in the Grothendieck group $K_0({\rm add}T)$, which transforms
in the same way as the index with respect to cluster-tilting objects. The main result (Theorem \ref{15}) states that all tropical
friezes on ${\mathcal{C}}_Q$ with $Q$ a Dynkin quiver are of the form $f_{T,m}$.
A different approach of this fact is given in Section 5 of \cite{FG}.
As an application, we show that for any tropical frieze $f$
on ${\mathcal{C}}_Q$, there exist cluster-tilting objects $T'$ (resp. $T''$) such that $f$ simultaneously takes non-negative (resp. non-positive)
values on all its indecomposable direct summands.

Section 6 gives a proof of a conjecture of Ringel
for the universal form of cluster-additive functions $f$ on stable translation quivers ${\mathbb{Z}} \Delta$ with $\Delta$ a
simply laced Dynkin diagram,
namely, $f$ is a non-negative linear combination of pairwise `compatible' (in the sense of Ringel) cluster-hammock functions.

\subsection*{Acknowledgments} The author thanks Sophie Morier-Genoud for pointing out J. Propp's work \cite{JP}. She is supported by the China Scholarship Council (CSC).
This is part of her Ph.~D.~thesis under the supervision of Professor
Bernhard Keller. She is deeply grateful to him for his guidance,
patience and kindness. She also sincerely thanks Pierre-Guy
Plamondon for helpful discussions and Zhonghua Zhao for constant
encouragement.

\vspace{.3cm}

\section{first properties of tropical friezes}

In this section, we recall Coxeter-Conway's frieze patterns at the beginning, then inspired by a tropicalized version of Coxeter-Conway's
frieze patterns of integers, we introduce tropical friezes on 2-Calabi-Yau triangulated categories. Apart from studying their first properties,
we also investigate some links between tropical friezes and cluster characters.

\subsection{Frieze patterns}

In early 1970s, Coxeter and Conway studied frieze patterns and triangulated polygons in \cite{Cox,coco1,coco2}. A frieze pattern
${\mathcal{F}}_n$ of order $n$ consists of $n-1$ infinite rows of positive numbers, whose first and last rows are filled with $1$.
Besides, the essential point is the {\it unimodular rule}, that is, for every four adjacent numbers in ${\mathcal{F}}_n$ forming a
diamond shape
$$\xymatrix @-1.8pc {
 &b &  \\
 a&  &d \\
 &c & }$$
the relation $ad = bc + 1$ is satisfied.
For example, the following diagram is a frieze pattern of order 6:
$$\xymatrix @-2pc{
\ldots & 1 & & 1 & & 1& & 1 & & 1 & & 1 & & 1 & & 1 & \ldots\\
& \ldots  & 2 & & 2 & & 2 & & 1& & 4& & 1 & & 2 & & 2 & \ldots \\
& & \ldots & 3 & & 3 & & 1& & 3 & & 3 & & 1 & & 3 & & 3 &\ldots\\
& & & \ldots & 4 & & 1 & & 2 & & 2 & & 2 & & 1 & & 4 & & 1 & \ldots\\
& & & & \ldots & 1 & & 1 & & 1& & 1& & 1& & 1& & 1& & 1& \ldots
}$$
A notable property of ${\mathcal{F}}_n$ is its periodicity with period a divisor of $n$.
More precisely, it is invariant under a glide reflection $\sigma$ which is $[{\frac{n}{2}}]$ times horizontal translation composed with
a horizontal reflection.

A frieze pattern ${\mathcal{F}}_n$ is determined by the elements in one of its diagonals (say $b_1 =1$, $b_2, \ldots, b_{n-2},
b_{n-1} = 1$), and it consists of integers if and only if $b_s$ divides $b_{s-1} + b_{s+1}$ for $s = 2, \ldots, n-2$. Let $a_0 = b_2$
and $a_1, a_2, \ldots$ be the numbers lying to the right of $a_0$ in the second row. Then we have
$$a_s = \frac{b_s + b_{s+2}}{b_{s+1}}, \qquad 1 \leq s \leq n-3.$$ A frieze pattern ${\mathcal{F}}_n$ can also be derived from $a_0,
\ldots, a_{n-4}$, since $a_{n-3}$ satisfies the linear equation
\begin{displaymath}
\left| \begin{array}{cccccc}
a_0 & 1& 0 & \ldots & 0 & 0  \\
1& a_1& 0 & \ldots& 0 & 0 \\
0 & 1 & a_2 & \ldots & 0 & 0\\
 &\ldots & & \ldots & & \\
0 & 0& 0& \ldots & 1 & a_{n-3}
\end{array} \right| = b_{n-1} = 1
\end{displaymath}
and ${\mathcal{F}}_n$ is symmetrical by the glide reflection $\sigma$. Moreover, ${\mathcal{F}}_n$ consists of integers if and only if $a_0,
\ldots, a_{n-4}$ are integers.

Let ${\mathcal{P}}_n$ be a regular $n$-gon with vertices $0, \ldots, n-1$. A triangulation $T$ of ${\mathcal{P}}_n$ is a maximal
set of non-crossing diagonals of ${\mathcal{P}}_n$, whose cardinality is always equal to $n-3$. Such a pair $({\mathcal{P}}_n, T)$
is called a triangulated $n$-gon. Let $a_r$ denote the number of triangles at vertex $r$ with respect to some triangulation $T$. Then
$$\ldots \quad  a_0 \quad a_1 \,  \ldots \, a_{n-1} \quad a_0\quad  \ldots$$ is the second row of a frieze pattern of integers.
Furthermore, the frieze patterns of integers of order $n$ are in bijection with triangulated $n$-polygons. %Let $T$ be a triangulation
%of ${\mathcal{P}}_n$ such that the associated quiver $Q_T$ defined in \cite{CCS} is an orientation of Dynkin type $A_{n-3}$. The cluster
%variables in the cluster algebra ${\mathcal{A}}_{Q_T}$ are parameterized by the diagonals of ${\mathcal{P}}_n$. When the initial
%cluster variables $u_1, \ldots, u_{n-3}$ (parameterized by the diagonals in $T$) are specialized to $1$, the frieze pattern
%of integers corresponding to $({\mathcal{P}}_n, T)$ is shown \cite{CC06} to be obtained from the Auslander-Reiten quiver of the cluster
%category ${\mathcal{C}}_{Q_T}$ by replacing each indecomposable object $M$ with
%$${\mathbf{x}}_M = CC(M)|_{u_1 = \ldots = u_{n-3} = 1}$$ where $CC(M)$ is the cluster variable of ${\mathcal{A}}_{Q_T}$ given by
%the image of $M$ under the Caldero-Chapoton map
%$$CC: obj({\mathcal{C}}_Q) \longrightarrow {\mathbb{Q}}(x_1, \ldots, x_n)$$ defined in \cite{CC06}.
%For example, the above given frieze pattern of order 6 corresponds to the following triangulated 6-gon
%$$\xymatrix@-0.5pc{
%& \ar@{-}[r] \ar@{.}[dd] & \ar@{-}[dr] &\\
% \ar@{-}[ur] & & & \ar@{-}[dl] \\
%& \ar@{-}[ul] \ar@{-}[r] \ar@{.}[uur] \ar@{.}[urr] & &
%}$$

Associated with an acyclic quiver $Q$, the authors observed in \cite{CC06} a generalized version of Coxeter-Conway's frieze patterns.
The elements of the generalized frieze pattern ${\mathcal{F}}_Q$ associated with $Q$ are cluster variables in the cluster algebra ${\mathcal{A}}_Q$.
Moreover, the sequences in ${\mathcal{F}}_Q$ satisfy linear recurrence relations if and only if
$Q$ is of Dynkin or affine type (see \cite{FZII,ARS,KS}). Of course,
there are more connections between frieze patterns and cluster algebras (see for instance \cite{FM,AD,BM}).

%A $R$-frieze $f$ on a stable translation quiver $(\Gamma, \tau)$ is a map from the vertex set $\Gamma_0$ to a commutative ring $R$,
%which satisfies the {\it mesh} relation (a generalization of the unimodular rule)
%$$f(\tau(i)) f(i) = \prod_{ (a: j \rightarrow i )\in {\Gamma}_1} f(j) +1.$$

The {\it{tropical semifield}} $({\mathbb{Z}}, \odot, \oplus)$ is the set $\mathbb{Z}$ of integers with multiplication and addition given by
$$a \odot b = a+b, \quad \quad a \oplus b = \max \{a, b \}.$$ Notice that the unit in the tropical semifield with respect to
the given multiplication is the number 0.

If we view the unimodular rule as an equation in the tropical semifield, then it becomes
$$a+d = \max \{b+c, 0\},$$
which is deduced from $$a \odot d = a+d \quad \mbox{and} \quad (b \odot c) \oplus 1 = \max \{b+c, 0\}.$$

\begin{example}
One can easily check that for every adjacent numbers $a, b, c, d$
forming a diamond shape with $a$ left and $d$ right in the following diagram
$$\xymatrix @-2.1pc{
\ldots & 0 & & 0 & & 0& & 0 & & 0 & & 0 & & 0 & & 0& \ldots\\
& \ldots  & 2 & & 1& & 1 & & -1& & 4& & -2 & & 2 & & 1 & \ldots \\
& & \ldots & 3 & & 2 & & -2& & 3 & & 2 & & -2 & & 3 & & 2 &\ldots\\
& & & \ldots & 4 & & -2 & & 2 & & 1 & & 1 & & -1 & & 4 & & -2 & \ldots\\
& & & & \ldots & 0 & & 0 & & 0& & 0& & 0& & 0& & 0& & 0& \ldots }$$
the relation $a+d = \max \{b+c, 0\}$ is satisfied. Notice that if we
omit the first and last rows which are filled with $0$, nothing will
change. We call such a diagram a {\em tropicalized frieze pattern}
of order 6. This diagram is also periodic with period a divisor of
6, it is also invariant under the same glide reflection $\sigma$ (as
frieze patterns). In fact, this is a general phenomenon: every
tropicalized frieze pattern of order $n$ is periodic. We will
explain this fact right after Proposition \ref{16}.
\end{example}

In the following, we will study tropical friezes on 2-Calabi-Yau
triangulated categories, especially on the cluster categories
associated with Dynkin diagrams. As we will see after Proposition
\ref{16}, this generalizes the above tropicalization of frieze
patterns of integers.

\subsection{Definitions and first properties}

Let $k$ be an algebraically closed field. Let $\mathcal{C}$ be a $k$-linear triangulated category
with suspension functor $\Sigma$ where all idempotents split. We further assume that the category $\mathcal{C}$

\begin{itemize}

\item[a)] is Hom-finite, {\it i.e.} the morphism space ${\mathcal{C}}(X,Y)$ is finite-dimensional
for any two objects $X$, $Y$ in $\mathcal{C}$ (which implies that $\mathcal{C}$ is Krull-Schmidt);

\item[b)] is 2-Calabi-Yau, {\it i.e.} there exist bifunctorial isomorphisms
$$D{\mathcal{C}}(X,Y) \simeq {\mathcal{C}}(Y,{\Sigma}^2X),  \quad X,\,Y \in {\mathcal{C}},$$
where $D$ denotes the duality functor ${\rm Hom}_k(?,k)$;

\item[c)] admits a cluster-tilting object $T$, {\it i.e.}
\begin{itemize}
\item[i)] $T$ is rigid (that is, ${\mathcal{C}}(T,\Sigma T) = 0$), and $T$ is basic
(that is, its indecomposable direct summands are
pairwise non-isomorphic),
\item[ii)] for each object $X$ of $\mathcal{C}$, if ${\mathcal{C}}(T,\Sigma X)$ vanishes, then $X$ belongs to the
subcategory add$T$ of direct summands of finite direct sums of copies of $T$.
 \end{itemize}

\end{itemize}

If a category ${\mathcal{C}}$ satisfies all these assumptions,
we say that ${\mathcal{C}}$ is a 2{\it{-Calabi-Yau category with
cluster-tilting object}}. A typical class of such categories is the class of
cluster categories \cite{BMRRT06} of connected finite
acyclic quivers. Throughout this article, our category $\mathcal{C}$ is
always a 2-Calabi-Yau category with cluster-tilting object.

\begin{defn}\label{1}
 A tropical frieze on $\mathcal{C}$ with values in the integer ring $\mathbb{Z}$ is a map
$$f: \, obj({\mathcal{C}}) \rightarrow {\mathbb{Z}}$$
such that
\begin{itemize}
\item[d1)] $f(X) = f(Y)$ if $X$ and $Y$ are isomorphic,
\item[d2)] $f(X \oplus Y) = f(X) + f(Y)$ for all objects $X$ and $Y$,
\item[d3)] for all objects $L$ and $M$ such that ${\rm dim Ext}^1_{\mathcal{C}} (L,M) = 1$, the equality
$$f(L) + f(M) = \max \{f(E), f(E') \}$$ holds, where $E$ and $E'$ are the middle terms of the non-split triangles
$$L \rightarrow E \rightarrow M \rightarrow \Sigma L \quad \mbox{and} \quad
 M\rightarrow E' \rightarrow L \rightarrow \Sigma M$$ with end terms $L$ and $M$.
\end{itemize}

Let $f$ and $g$ be two tropical friezes on the same category $\mathcal{C}$. The sum $f+g$ clearly satisfies items d1) and d2).
For item d3), we have that
\begin{align*}
(f+g)(L) + (f+g)(M) & = (f(L) + f(M)) + (g(L) + g(M)) \\
       & = \max \{f(E), f(E') \} + \max \{g(E), g(E') \}.
\end{align*}
Then $f+g$ is a tropical frieze if and only if for all pairs $(E, E')$ as in item d3) the equality
$$\max \{f(E), f(E') \} + \max \{g(E), g(E') \} = \max \{(f+g)(E), (f+g)(E') \}$$ holds.
Notice that for two integers $a,\,b$, the number $$\max \{a, b\} = \frac{a+b + |a-b|}{2}.$$
Thus, the sum $f+g$ is a tropical frieze if and only if for all pairs $(E, E')$ as in item d3) the equality
$$|f(E) - f(E')| + |g(E) - g(E')| = |(f(E) - f(E')) + (g(E) - g(E'))| $$ holds,
if and only if the inequality
$$(f(E) - f(E'))(g(E) - g(E')) \geq 0$$ holds. If two tropical friezes satisfy such a property, then we say that
they are {\it compatible}.

Now we state a simple property of tropical friezes.

\begin{prop}
Let $f_1, \ldots, f_n$ be tropical friezes on the same category $\mathcal{C}$. Then the sum ${\sum}_i f_i$ is a tropical frieze
if and only if the functions $f_i$ are pairwise compatible.
\end{prop}

\begin{pf}
This statement is a trivial generalization of the above analysis:

the sum ${\sum}_i f_i$ is a tropical frieze if and only if for all pairs $(E, E')$ as in item d3) the equality
$$ {\sum}_i |f_i(E) - f_i(E')| = | {\sum}_i (f_i(E) - f_i(E'))|$$ holds, if and only if $f_i(E) - f_i(E')$
are simultaneously non-negative or simultaneously non-positive for all integers $1 \leq i \leq n$, if and only if the tropical friezes $f_i$
are pairwise compatible.
\end{pf}

Similarly, one can obtain that the difference $f-g$ is a tropical frieze if and only if
for all pairs $(E, E')$ as in item d3) the equality
$$|f(E) - f(E')| - |g(E) - g(E')| = |(f(E) - f(E')) - (g(E) - g(E'))|$$ holds, if and only of the inequalities
$$|f(E) - f(E')| \geq |g(E) - g(E')| \quad \mbox{and} \quad
(f(E) - f(E'))(g(E) - g(E')) \geq 0$$ hold. If two tropical friezes satisfy such a property, then we say that
they are {\it strongly compatible}.

Let ${\mathcal{C}}_Q$ be the cluster category of a Dynkin quiver $Q$.
For any indecomposable object $X$ of ${\mathcal{C}}_Q$, the space ${\rm Hom}_{{\mathcal{C}}_Q}(X,X)$ is one-dimensional,
so we have that ${\rm dim}{\rm Ext}^1_{{\mathcal{C}}_Q}(\Sigma X, X) = 1$. The associated non-split triangles are of the following form
$$\qquad \Sigma X \rightarrow E \stackrel{g}\rightarrow X \rightarrow \Sigma^2 X \quad \mbox {and}
\quad X \rightarrow 0 \stackrel{g'}\rightarrow \Sigma X \rightarrow \Sigma X.    \quad \qquad  (\ast)$$

The following proposition is quite similar to the statements for cluster-additive functions
on stable translation quivers given in Section 1 of \cite{Rin11}.

\begin{prop}\label{4}
Let $Q$ be a Dynkin quiver. Then
any tropical frieze on ${\mathcal{C}}_Q$ which takes only non-positive values or only non-negative values is the zero function.
\end{prop}

\begin{pf}
Let $f$ be a non-zero tropical frieze on ${\mathcal{C}}_Q$ with non-positive values and
$X$ an indecomposable object such that $f(X) < 0$. From the non-split triangles ($\ast$) above, we have that
$$f(\Sigma X) = \max \{ f(E), 0\} - f(X) \geq 0 - f(X) > 0,$$ which is a contradiction.
Therefore, any tropical frieze with only non-positive values is the zero function.

Let $f$ be a tropical frieze on ${\mathcal{C}}_Q$ with non-negative values. We lift $f$ in the natural way to a
function $f'$ which is $(\tau^{-1} \Sigma)$-invariant on the bounded derived category ${\mathcal{D}}_Q$ of the category ${\rm mod}kQ$. Here $\tau$
is the Auslander-Reiten translation on ${\mathcal{D}}_Q$.
Denote by $\phi$
the canonical equivalence \cite{Hap87} from the mesh category
of the translation quiver ${\mathbb{Z}}Q$ to the full subcategory ${\rm ind}({\mathcal{D}}_Q)$ of indecomposables of ${\mathcal{D}}_Q$.
We define a function $f''$ on ${\mathbb{Z}}Q$ by setting $f''= f'\phi$.
Let $z$ be any vertex of ${\mathbb{Z}}Q$. In ${\mathcal{D}}_Q$ we have the Auslander-Reiten triangle \cite{Hap} as follows
$$\phi(\tau z) \rightarrow \bigoplus_{y\rightarrow z} \phi(y) \rightarrow \phi(z) \rightarrow \Sigma \phi(\tau z),$$
where `$y \rightarrow z$' in the middle term are arrows in ${\mathbb{Z}}Q$.
Its image (still use the same notation) in ${\mathcal{C}}_Q$ is a non-split triangle. The other non-split triangle
with end terms $\phi({z})$ and $\phi(\tau z)$ in ${\mathcal{C}}_Q$ is
$$\phi(z) \rightarrow 0 \rightarrow \phi(\tau z) \stackrel{\simeq}\longrightarrow \Sigma \phi(z).$$ Hence, we can deduce that
\begin{align*}
 f''(\tau z) + f''(z) & = f'(\phi(\tau z)) + f'(\phi(z)) = \max \{ \sum_{y\rightarrow z} f'(\phi(y)),0 \} \\
 &  = \sum_{y \rightarrow z} f'(\phi(y)) = \sum_{y \rightarrow z} f''(y).
\end{align*}
As a consequence, the function $f''$ is an additive function on ${\mathbb{Z}}Q$ with non-negative values, which implies that $f''$
is the zero function \cite{HPR}. Therefore, the function $f$ is the zero function on ${\mathcal{C}}_Q$.
\end{pf}

\subsection{Cluster characters and tropical friezes}

In this subsection, we will see some links between cluster characters and tropical friezes.

Let d$2'$) denote the item obtained from item d2) in Definition \ref{1} in which the equality becomes $f(X \oplus Y) = f(X) f(Y)$, and
d$3'$) the item obtained from item d3) in Definition \ref{1} in which the equality becomes $f(L) f(M) = f(E) + f(E')$.
A map $\chi: obj({\mathcal{C}}) \rightarrow A$, where $A$ is a commutative ring,
is called a {\it{cluster character}} in \cite{Pal08} if it satisfies items d1), d$2'$) and d$3'$).

\begin{rem} \label{13}
Let $\chi$ be a cluster character mapping from ${\mathcal{C}}$ to the tropical semifield $({\mathbb{Z}}, \odot, \oplus)$. Then we obtain
the following equalities
$$\chi (X \oplus Y) = \chi(X) \odot \chi(Y) = \chi(X) + \chi(Y),$$
$$\chi(L) + \chi(M) = \chi(L) \odot \chi(M) = \chi(E) \oplus \chi(E') = \max \{ \chi(E), \chi(E') \}. $$
As a result, the map $\chi$ is a tropical frieze mapping to the integer ring ${\mathbb{Z}}$.
\end{rem}

Let $Q$ be a connected finite acyclic quiver with vertex set $\{ 1, \ldots, n \}$ and ${\mathcal{C}}_Q$ its associated cluster category. It was proved, in \cite{CK08} for Dynkin quivers
and in \cite{CK06} for acyclic quivers, that the Caldero-Chapoton map
$$CC: obj({\mathcal{C}}_Q) \longrightarrow {\mathbb{Q}}(x_1, \ldots, x_n)$$ defined in \cite{CC06} is a
cluster character.

\begin{example} \label{12}
Let $X$ be an object of ${\mathcal{C}}_Q$. Then the image $CC(X)$ can be written uniquely as
$$CC(X) = \frac{ h(x_1, \ldots, x_n)}{\prod_{i=1}^n x_i^{d_i(X)} },$$ where the polynomial $h(x_1, \ldots, x_n)$
is not divisible by any $x_i, \, 1 \leq i \leq n$. Look at the
function $$d_i: obj({\mathcal{C}}_Q) \rightarrow {\mathbb{Z}}$$ with
$d_i(X)$ given as in the above expression for each object $X$ of
${\mathcal{C}}_Q$.

We use elementary properties of polynomials.
From the equlity $CC(X \oplus Y) = CC(X) CC(Y)$
in item d$2'$), one clearly sees that $d_i (X \oplus Y) = d_i(X) + d_i(Y)$.
It is also not hard to calculate
the denominators of the two hand sides of the equality $CC(L)CC(M) = CC(E)+CC(E')$ in item d$3'$), which gives us the equality
$d_i(L) + d_i(M) = \max \{d_i(E), d_i(E')\}$. Therefore, each function $d_i$ is a tropical frieze on ${\mathcal{C}}_Q$.
\end{example}

However, the sum $d_i + d_j$ is not always a tropical frieze on ${\mathcal{C}}_Q$. We choose the linear orientation of $A_3$.
The Auslander-Reiten quiver of the cluster category ${\mathcal{C}}_{\vec{A_3}}$ is
$$
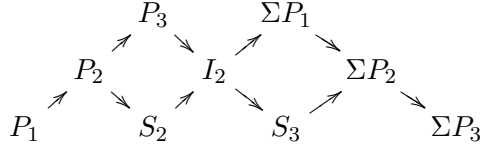
\xymatrix @-1.5pc {
 & & P_3 \ar[dr] & & \Sigma P_1 \ar[dr] & &  \\
 & P_2 \ar[ur] \ar[dr] & & I_2 \ar[ur] \ar[dr] & & \Sigma P_2 \ar[dr] & \\
 P_1 \ar[ur] & & S_2 \ar[ur] & & S_3 \ar[ur] & & \Sigma P_3  }$$ where $P_i$ (resp. $I_i$, $S_i$)
is the right projective (resp. injective, simple) $k \vec{A_3}$-module associated to vertex $i$. By definition $CC(\Sigma P_2) = x_2$
and one can calculate that
$$CC(P_1) = \frac{1+x_2}{x_1}, \,\, CC(S_3) = \frac{1+x_2}{x_3}, \,\, CC(P_3) = \frac{x_1 + x_1x_2 + x_3 +x_2x_3}{x_1x_2x_3}.$$
The space ${\rm Ext}^1_{{\mathcal{C}}_{\vec{A_3}}}(\Sigma P_2, P_3)$ is 1-dimensional and the non-split triangles are
$$\Sigma P_2 \rightarrow P_1 \rightarrow P_3 \rightarrow I_2 \quad \mbox{and} \quad P_3 \rightarrow S_3 \rightarrow \Sigma P_2 \rightarrow \Sigma P_3.$$
Consider the function $d_1 + d_3$. We have that
\begin{align*}
 (d_1 + d_3)(\Sigma P_2) + (d_1 + d_3)(P_3) & = (0+0)+(1+1) = 2 \\
 \max \{(d_1 + d_3)(P_1), (d_1 + d_3)(S_3) \}& = \max \{1+0, 0+1 \} = 1.
\end{align*}
Thus, the sum $d_1 +d_3$ is not a tropical frieze. In another way,
since $$(d_1(P_1) - d_1(S_3)) (d_3(P_1) - d_3(S_3)) = (1-0)(0-1) = -1 < 0, $$
the tropical friezes $d_1$ and $d_3$ are not compatible.
As a consequence, the difference $d_1 - d_3$ is not a tropical frieze on ${\mathcal{C}}_{\vec{A_3}}$ either.
\smallskip

Let $T$ be a cluster-tilting object of ${\mathcal{C}}$ and $T_1$ an
indecomposable direct summand of $T$. Iyama and Yoshino proved in
\cite{IY08} that, up to isomorphism, there is a unique
indecomposable object $T_1^{\ast}$ not isomorphic to $T_1$ such that
the object $\mu_1 (T)$ obtained from $T$ by replacing the
indecomposable direct summand $T_1$ with $T_1^{\ast}$ is
cluster-tilting. We call $\mu_1(T)$ the {\it{mutation}} of $T$ at
$T_1$. There are non-split triangles, unique up to isomorphism,
$$T_1^{\ast} \rightarrow E \rightarrow T_1 \rightarrow \Sigma T_1^{\ast} \quad {\mbox{and}} \quad
T_1 \rightarrow E' \rightarrow T_1^{\ast} \rightarrow \Sigma T_1$$
with $E$ and $E'$ in add$(T/T_1)$. A category $\mathcal{C}$ is said
to be {\it{cluster-transitive}} if any two basic cluster-tilting
objects of ${\mathcal{C}}$ can be obtained from each other by a
finite sequence of mutations.

The following property of tropical friezes on a cluster-transitive category ${\mathcal{C}}$
is quite similar to that \cite{Pal08} of cluster
characters on ${\mathcal{C}}$.

\begin{prop} \label{11}
Let ${\mathcal{C}}$ be a cluster-transitive category and $T = T_1 \oplus \ldots \oplus T_n$ a basic cluster-tilting object of $\mathcal{C}$.
Suppose that $f$ and $g$ are two tropical friezes on $\mathcal{C}$ such that $f(T_i) = g(T_i), \, 1 \leq i \leq n$. Then $f$ and $g$
coincide on all subcategories add$T'$, where $T'$ is any cluster-tilting object of $\mathcal{C}$.
\end{prop}

\begin{pf}

By assumption we know that $f$ and $g$ coincide
on all indecomposable direct summands of $T$.
We will prove this proposition by recursion on the minimal number of
mutations linking a basic cluster-tilting object to $T$.

Now let $T' = T'_1 \oplus \ldots \oplus T'_n$ be a basic cluster-tilting object satisfying
that $f(T'_i) = g(T'_ i)$ for all integers $1 \leq i \leq n.$ Assume that
$T{''} = \mu_1(T') = T^{''}_1 \oplus T'_2 \oplus \ldots \oplus T'_n$ is the mutation of $T'$ in
direction 1. Then we have the non-split triangles
$$T^{''}_1 \rightarrow E \rightarrow T'_1 \rightarrow \Sigma T^{''}_1 \quad \mbox{and} \quad
 T'_1\rightarrow E' \rightarrow T^{''}_1 \rightarrow \Sigma T'_1$$ with middle terms $E$ and $E'$ both belonging to add($T' / T'_1$).
Hence, the following equlities
$$f(T^{''}_1) = \max \{f(E), f(E') \} - f(T'_1) = \max \{g(E), g(E') \} - g(T'_1) = g(T^{''}_1)$$ hold.
This completes the proof.
\end{pf}

Let ${\mathcal{C}}_Q$ be the cluster category of a connected finite acyclic quiver $Q$.
It was shown in \cite{BMRRT06} that ${\mathcal{C}}_Q$ is cluster-transitive and any rigid indecomposable object
of ${\mathcal{C}}_Q$ is a direct summand of a cluster-tilting object. If $f$ and $g$ are two tropical friezes
on ${\mathcal{C}}_Q$ which coincide on
all indecomposable direct summands of some cluster-tilting object, by Proposition \ref{11},
they coincide on all rigid objects. In particular, when $Q$ is Dynkin, the two tropical friezes $f$ and $g$ are equal.

\vspace{.3cm}

\section{tropical friezes from indices}

\subsection{Reminder on indices}
Let $X$ be an object of $\mathcal{C}$ and $T$ a cluster-tilting object of $\mathcal{C}$. Following \cite{KR07}, we have triangles
$$T^{X}_1 \rightarrow T^{X}_0 \rightarrow X \rightarrow \Sigma T^{X}_1 \quad \mbox{and} \quad
X \rightarrow \Sigma^2 T^0_X \rightarrow \Sigma^2 T^1_X  \rightarrow \Sigma X,$$where $T^X_1,\,T^X_0,\,T^0_X$ and $T^1_X$ belong to add$T$.
Recall that the index and coindex of $X$ with respect to $T$ are defined to be the classes in the split Grothendieck group $K_0$(add$T$) of
the additive category add$T$ as follows
$${\rm ind}_T(X) = [T^X_0] - [T^X_1] \quad \mbox{and} \quad {\rm coind}_T(X) = [T^0_X]
- [T^1_X] \, ,$$ which do not depend on the choices of the above
triangles.

Assume that $T$ is the direct sum of $n$ pairwise non-isomorphic indecomposable objects $T_1,\,\ldots, T_n$.
Let $B$ be the endomorphism algebra of $T$ over ${\mathcal{C}}$. We denote the indecomposable right projective
$B$-module ${\mathcal{C}}(T,T_i)$ by $P_i$ and denote its simple top by $S_i$. For any two finite-dimensional $B$-modules $X$ and $Y$, set
\begin{align*}
 \langle X,Y \rangle & = {\rm dim} \,{\rm Hom}_B (X,Y) - {\rm dim} \,{\rm Ext}^1_{B}(X,Y) \quad \mbox{and}\\
\langle X,Y \rangle_a & = \langle X,Y \rangle - \langle Y,X \rangle.
\end{align*}
In \cite{Pal08} Palu has proved that $\langle \, , \, \rangle_a$
is a well-defined bilinear form on the Grothendieck group $K_0 ({\rm mod} B)$ of the abelian category ${\rm mod}B$ of finite-dimensional right $B$-modules.
Let $F$ denote the functor ${\mathcal{C}}(T,?)$. It was shown in \cite{KR07} that $F$ induces an equivalence of categories
$${\mathcal{C}}/{\rm add}(\Sigma T) \stackrel{\simeq}\longrightarrow {\rm mod}B.$$
%where $(\Sigma T)$ is the ideal of morphisms of $\mathcal{C}$ which factor through an object of add($\Sigma T$).
Let $m$ be a class in $K_0 ({\rm mod}B)$. We define a function $f_{T,m}$ from $\mathcal{C}$ to ${\mathbb{Z}}$ as
$$f_{T,m}(X) = \langle F({\rm ind}_T(X)), m \rangle, \quad X \in {\mathcal{C}}.$$ When it does not cause confusion, we simply write ${\rm ind}(X)$
instead of ${\rm ind}_T(X)$.

\subsection{Tropical friezes}

In this subsection, we will give a sufficient condition for the function $f_{T,m}$ to be a tropical frieze on $\mathcal{C}$.
Moreover, when ${\mathcal{C}} = {\mathcal{C}}_Q$ the cluster category of a Dynkin quiver $Q$, we will see that this sufficient condition is also a necessary condition.

\begin{thm}\label{2}
Assume that $\langle S_i, m \rangle_a \geq 0$ for each simple $B$-module $S_i \, (1 \leq i \leq n)$. Then the function $f_{T,m}$ is a tropical frieze.
\end{thm}

\begin{pf}
The function $f_{T,m}$ clearly satisfies the terms d1) and d2) in Definition \ref{1}. Now Let $L$ and $M$ be objects of $\mathcal{C}$
such that ${\rm dim Ext}^1_{\mathcal{C}} (L,M) = 1$. Let
$$L \stackrel{h}\rightarrow E \stackrel{g}\rightarrow M \rightarrow \Sigma L \quad \mbox{and} \quad
 M \stackrel{h'}\rightarrow E' \stackrel{g'}\rightarrow L \rightarrow \Sigma M$$ be the associated non-split triangles.

First, let $C \in {\mathcal{C}}$ be any lift of ${\rm Coker}(Fg)$. We know from \cite{Pal08} that
\begin{align*}
 {\rm ind}(E) & = {\rm ind}(L) + {\rm ind}(M) - {\rm ind}(C) - {\rm ind}(\Sigma^{-1}C) \quad {\mbox{and}} \\
\langle FC, m \rangle_a & = \langle F({\rm ind}(C)), m \rangle + \langle F({\rm ind}(\Sigma^{-1}C)), m \rangle.
\end{align*}
By assumption $\langle S_i, m \rangle_a \geq 0$ for each simple $B$-module $S_i \, (1 \leq i \leq n)$.
So we have that $\langle FC, m \rangle_a \geq 0$. Thus,
\begin{align*}
 \langle F({\rm ind}(E)), m \rangle & = \langle F({\rm ind}(L)), m \rangle + \langle F({\rm ind}(M)), m \rangle - \langle FC, m \rangle_a \\
 & \leq \langle F({\rm ind}(L)), m \rangle + \langle F({\rm ind}(M)), m \rangle.
\end{align*}
Similarly, we obtain another inequality
$$\langle F({\rm ind}(E')), m \rangle \leq \langle F({\rm ind}(L)), m \rangle + \langle F({\rm ind}(M)), m \rangle.$$
It follows that $$\max \{ f_{T,m}(E), f_{T,m}(E') \} \leq f_{T,m}(L) + f_{T,m}(M).$$

Second, we consider the identity maps $id_M: M \rightarrow M$ and $id_L: L \rightarrow L$.
Thanks to the dichotomy phenomenon shown in \cite{Pal08}, exactly one
of the conditions $FM = (Fg)(F E)$ and $FL = (Fg')(F E')$ is true. Assume that the first condition holds, then $Fg$ is an epimorphism
and $FC$ vanishes. %which means that ${\rm Coker(Fg')}$ is non-zero.
Therefore, we have that
\begin{align*}
f_{T,m}(E) & = \langle F({\rm ind}(E)), m \rangle = \langle F({\rm ind}(L)), m \rangle + \langle F({\rm ind}(M)), m \rangle \\
& = f_{T,m}(L) + f_{T,m}(M).
\end{align*}
As a consequence, the equality $$f_{T,m}(L) + f_{T,m}(M) =
\max \{ f_{T,m}(E), f_{T,m}(E') \}$$ holds and $f_{T,m}$ is a tropical frieze.
\end{pf}

\subsection{Another proof}

For $L \in {\mathcal{C}}$ and $e \in {\mathbb{N}}^n$, we denote by $Gr_e({\rm Ext}^1_{\mathcal{C}}(T,L))$ the quiver Grassmannian of $B$-submodules
of the $B$-module ${\rm Ext}^1_{\mathcal{C}}(T,L)$ whose dimension vector is $e$ and by $\chi(Gr_e({\rm Ext}^1_{\mathcal{C}}(T,L)))$ its Euler-Poincar${\rm \acute{e}}$ characteristic
for ${\rm \acute{e}}$tale cohomology with proper support.

For $1 \leq i \leq n$, we define the integer $g_i (L)$ to be the multiplicity of $[T_i]$ in the index ${\rm ind}(L)$ and define
the element $X'_L$ of the field ${\mathbb{Q}}(x_1, \ldots, x_n)$ by
$$X'_L = \prod^n_{i = 1} x_i^{g_i(L)} \sum_{e} {\chi}(Gr_e({\rm Ext}^1_{\mathcal{C}}(T,L)))\prod^n_{i=1}x_i^{\langle S_i,e \rangle_a},$$
where the sum ranges over all tuples $e \in {\mathbb{N}}^n$. This is
a vastly generalized form of the $CC$ map. It was proved in
\cite{Pal08} that the function $X'_?$ is a cluster character from
$\mathcal{C}$ to ${\mathbb{Q}}(x_1, \ldots, x_n)$. If we define
functions $d_i$ on ${\mathcal{C}}$ as in Example \ref{12} by
replacing $CC$ map with $X'_?$, then each function $d_i$ is also a
tropical frieze.
%Let $d_i (L)$ be the power of $x_i$ such that $$X'_L = \prod_{i=1}^n x_i^{d_i(L)} \cdot f(x_1, \ldots, x_n),$$
%where the polynomial $f(x_1, \ldots, x_n)$
%is not divisible by any $x_i, \, 1 \leq i \leq n$. Similarly as in Example \ref{12}, the function $d_i$ is a tropical frieze
%from $\mathcal{C}$ to the integer ring $\mathbb{Z}$.

We will use the tropical semifield (${\mathbb{Z}}, \odot, \oplus$)
to give another proof of Theorem \ref{2} for ${\mathcal{C}} =
{\mathcal{C}}_Q$ where $Q$ is a Dynkin quiver with $n$ vertices. Let
$T$ be a cluster-tilting object of $\mathcal{C}$ and $B$ its
endomorphism algebra. Notice that any indecomposable object of
$\mathcal{C}$ is a direct summand of some cluster-tilting object
which is obtained from $T$ by a finite sequence of mutations. Since
$X'_{T_i} = x_i$ and $X'_?$ is a cluster character, the image $X'_L$
lies in the universal semifield ${\mathbb{Q}}_{sf}(x_1, \ldots,
x_n)$ (section 2.1 in \cite{BFZ95}). For an element $m \in K_0({\rm
mod}B)$, we define the map $$\varphi_m: {\mathbb{Q}}_{sf}(x_1,
\ldots, x_n)  \longrightarrow ({\mathbb{Z}}, \odot, \oplus)$$ as the
unique homomorphism between semifields which takes $x_i = X'_{T_i}$
to the integer $\langle F({\rm ind} (T_i)), m \rangle$. Then the
composition ${\varphi}_m X'_?$ is a cluster character from
$\mathcal{C}$ to $({\mathbb{Z}}, \odot, \oplus)$ and thus a tropical
frieze from $\mathcal{C}$ to the integer ring $\mathbb{Z}$ by Remark
\ref{13}. When ${\mathcal{C}} = {\mathcal{C}}_Q$ with $Q$ a Dynkin
quiver, Nakajima \cite{N} showed that ${\chi}(Gr_e({\rm
Ext}^1_{\mathcal{C}}(T,L)))$ is a non-negative integer. Now we write
down explicitly the function
\begin{align*}
{\varphi}_m X'_L & = \max_{e} \{ \sum^n_{i=1} (g_i(L) + \langle S_i, e \rangle_a) \langle F({\rm ind}(T_i)), m \rangle \} \\
&  = \max_e \{ \langle F({\rm ind}(L)), m \rangle + \sum^n_{i=1} \langle S_i, e \rangle_a \langle F({\rm ind}(T_i)), m \rangle \} \\
& = \max_e \{ \langle F({\rm ind}(L)), m \rangle - \langle (g_i (\Sigma^{-1} Y) + g_i(Y)) \langle F({\rm ind}(T_i)), m \rangle \} \\
& = \max_e \{ \langle F({\rm ind}(L)), m \rangle - \langle F({\rm ind}(\Sigma^{-1}Y)+{\rm ind}(Y)), m \rangle \} \\
& = \max_e \{ \langle F({\rm ind}(L)), m \rangle - \langle FY, m \rangle_a \} \\
& = \max_e \{ \langle F({\rm ind}(L)), m \rangle - \sum^n_{i=1} e_i \langle S_i, m \rangle_a\}
\end{align*}
where $e$ ranges over all elements in $K_0 ({\rm mod}B)$
such that $\chi(Gr_e({\rm Ext}^1_{\mathcal{C}}(T,L)))$ is non zero and $Y$ is an object of $\mathcal{C}$ satisfying $FY = e = (e_i)_i \in K_0({\rm mod}B)$.
If $\langle S_i, m \rangle_a \geq 0$ for each simple $B$-module $S_i$, then we have that
$${\varphi}_m X'_L = \langle F({\rm ind}(L)), m \rangle = f_{T,m}(L).$$
Thus, the function $f_{T,m}$ is equal to ${\varphi}_m X'_?$ and is a tropical frieze.

\begin{rem}\label{3}
Let ${\mathcal{C}}_Q$ be the cluster category associated to a Dynkin
quiver $Q$. Let $T$ be a cluster-tilting object of ${\mathcal{C}}_Q$
and $B$ its endomorphism algebra. Let $F$ be the functor ${\rm
Hom}_{{\mathcal{C}}_Q}(T,?)$. In fact, the sufficient condition for
a function $f_{T,m}$ to be a tropical frieze in Theorem \ref{2} is
also a necessary condition in this situation.

For any indecomposable object $X$ of ${\mathcal{C}}_Q$, %the space ${\rm Hom}_{{\mathcal{C}}_Q}(X,X)$ is one-dimensional,
%then we have that ${\rm dim}{\rm Ext}^1_{{\mathcal{C}}_Q}(\Sigma X, X) = 1$. Assume the associated non-split triangles are
%$$\Sigma X \rightarrow E \stackrel{g}\rightarrow X \rightarrow \Sigma^2 X \quad \mbox {and}
%\quad X \rightarrow 0 \stackrel{g'}\rightarrow \Sigma X \rightarrow \Sigma X.    \quad \quad (\ast)$$
look at the second triangle associated to $X$ in ($\ast$) before Proposition \ref{4}, whose image under $F$ is
$$ FX \rightarrow 0 \stackrel{Fg'}\longrightarrow F(\Sigma X) \stackrel{\simeq}\longrightarrow F(\Sigma X).$$
We have that ${\rm Coker}(Fg') = F(\Sigma X)$. If $X$ does not belong to add$T$, then ${\rm Coker(Fg')}$ is not zero which implies that
${\rm Coker}(Fg)$ vanishes by the dichotomy phenomenon.
Let $m$ be a class in $K_0({\rm mod}B)$. From the proof of Theorem \ref{2}, we know that
$$f_{T,m}(E) = f_{T,m}(\Sigma X) + f_{T,m}(X) = \langle F(\Sigma X),m \rangle_a.$$ Assume that $f_{T,m}$ is a tropical frieze.
Then it follows that
$$f_{T,m}(\Sigma X) + f_{T,m}(X) = \max \{ f_{T,m}(E), 0 \} \geq 0.$$ Thus, for every indecomposable object $X \notin$ add$T$,
the value $\langle F(\Sigma X),m \rangle_a$ is non-negative, particularly when $F(\Sigma X)$ is a simple $B$-module $S_i$.
\end{rem}

\begin{example}
Let $Q$ be an acyclic quiver and $j$ a sink of $Q$ (that is, no arrows of $Q$ start at $j$). Let $T$ be the image of $kQ$ %$ = e_1 kQ \oplus \ldots \oplus e_n kQ$
in ${\mathcal{C}}_Q$ under the canonical inclusion and $B$ its endomorphism algebra ${\rm End}_{{\mathcal{C}}_Q}(T)$.
%where $e_i$ is the idempotent associated to vertex $i$.
%Let $m$ be the simple (${\rm End}_{{\mathcal{C}}_Q} (T)$)-module $S_j$.
For each simple $B$-module $S_i$, we have that
$$\, \, \langle S_i, S_j \rangle_a = - {\rm dim} {\rm Ext}^1_B (S_i, S_j) + {\rm dim} {\rm Ext}^1_B (S_j, S_i) = {\rm dim} {\rm Ext}^1_B (S_j, S_i)$$
\begin{center}
= the number of arrows from $i$ to $j$ in $Q \, \geq 0. \, \quad $
\end{center}
As an application of Theorem \ref{2}, the function $f_{T,S_j}$ is a tropical frieze.

Similarly, if $j$ is a source of an acyclic quiver $Q$, that is, no arrows of $Q$ end at $j$, then $f_{T, -S_j}$ is a tropical frieze.
\end{example}

Using a similar method as in the second proof of Theorem \ref{2}, it is not hard to get the following proposition:

\begin{prop}\label{16}
Let ${\mathcal{C}}_Q$ be the cluster category of a Dynkin quiver $Q$ and $T = T_1 \oplus \ldots \oplus T_n$
a basic cluster-tilting object of ${\mathcal{C}}_Q$.
Then the map
$$  \Phi_T : \{\mbox{tropical\, friezes\,on}\,\, {\mathcal{C}}_Q\} \longrightarrow {\mathbb{Z}}^n $$
given by $\Phi_T(f) = (f(T_1), \ldots, f(T_n))$
is a bijection.
\end{prop}

\begin{pf}
For any fixed $n$-tuple $\underline{a} = (a_1, \ldots, a_n)$ in ${\mathbb{Z}}^n$, there is a unique homomorphism of semifields
$$\phi_{\underline{a}}: {\mathbb{Q}}_{sf}(x_1, \ldots, x_n) \longrightarrow ({\mathbb{Z}},\odot, \oplus)$$ such that
$\phi_{\underline{a}}(x_i) = a_i$. We denote the composition $\phi_{\underline{a}} X'_?$ by $f_{\underline{a}}$.
Then $f_{\underline{a}}$ is a tropical frieze on ${\mathcal{C}}_Q$ satisfying $f_{\underline{a}}(T_i) = a_i$. Therefore, the
map $\Phi_T$ is a surjection. The injectivity follows from Proposition \ref{11}. Hence, the map $\Phi_T$ is bijective.
\end{pf}

Now we give an explanation of the periodicity phenomenon which is
stated at the end of subsection 2.1. Let ${\mathcal{F}}_n^t$ be a
tropicalized frieze pattern of order $n (> 3)$. Let $Q$ be a quiver
of type $A_{n-3}$. Then ${\mathcal{F}}_n^t$ gives a function
(denoted by $f$) on the Auslander-Reiten quiver $\Gamma$ of
${\mathcal{D}}_Q$. Each subquiver ($y$ or $z$ may not appear)
$$
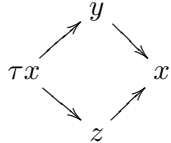
\xymatrix @-1pc {
 &y \ar[dr] &  \\
 \tau x \ar[ur] \ar[dr] &  &x \\
 &z \ar[ur]& }$$
in $\Gamma$ induces an Auslander-Reiten triangle in ${\mathcal{D}}_Q$
$$\tau x \rightarrow y \oplus z \rightarrow x \rightarrow \Sigma \tau x.$$
Since ${\mathcal{F}}_n^t$ is a tropicalized frieze pattern, the function $f$ satisfies that
$$f(\tau x) + f(x) = \max \{ f(y)+f(z), 0 \}.$$ %Therefore, for all non-split triangles as the triangles $(\ast)$
%before Proposition \ref{4},
%the function $f$ satisfies item d3) in Definition \ref{1}.

Let $\mathcal{S}$ be any slice in $\Gamma$. Set $T = \bigoplus_{y \in {\mathcal{S}}} y$.
Then the image of $T$ %under the canonical projection $\pi: {\mathcal{D}}_Q \rightarrow {\mathcal{C}}_Q$
is a basic cluster-tilting object of ${\mathcal{C}}_Q$.
By Proposition \ref{16}, there exists a unique tropical frieze $g: {\mathcal{C}}_Q \rightarrow {\mathbb{Z}}$ such that
$g(y) = f(y)$ for all $y \in {\mathcal{S}}$.
We extend $g$ in a natural way to a $(\tau^{-1} \Sigma)$-invariant function on ${\mathcal{D}}_Q$ (still denote as $g$). Then $g$
also satisfies the above equation as $f$. Therefore, the two functions $f$ and $g$ are equal. Moreover, for each integer $i$, we have that
$$f((\tau^{-1} \Sigma)^i x) = g ((\tau^{-1} \Sigma)^i x) = g(x) = f(x) \quad \mbox{and}$$
$$f((\tau^{-n})^i x) = g((\tau^{-n})^i x) = g((\tau^{-2} \tau^{2-n})^i x) = g((\tau^{-1} \Sigma)^{2i} x) = g(x).$$
In conclusion, ${\mathcal{F}}_n^t$ is periodic with period a divisor
of $n$, and it is invariant under the glide reflection $\sigma$.

\subsection{Cluster-hammock functions and tropical friezes}

In this subsection, we will see that the cluster-hammock functions defined by Ringel \cite{Rin11} always give rise to tropical friezes,
while their sums do not, even for pairwise `compatible' (in the sense of Ringel) cluster-hammock functions.

Let $\Gamma = {\mathbb{Z}} Q$ be the translation quiver of a Dynkin quiver $Q$. For any vertex $x$ of $\Gamma$,
Ringel \cite{Rin11} defined the {\em cluster-hammock function} $h_x: \Gamma_0 \rightarrow {\mathbb{Z}}$
by the following properties
\begin{itemize}
\item[a)] $h_x (x) = -1;$
\item[b)] $ h_x(y) = 0$ for $y \neq x \in {\mathcal{S}},$ where $\mathcal{S}$ is any slice containing $x$;
\item[c)] $h_x(z)+h_x(\tau z)=\sum_{y \rightarrow z} \max \{h_x(y), 0 \}$ for all $z \in \Gamma_0$.
\end{itemize}
As shown in \cite{Rin11}, the cluster-hammock function $h_x$ is $({\tau}^{-1} \Sigma)$-invariant
and takes the value $-1$ on the $({\tau}^{-1} \Sigma)$-orbit of $x$ while it takes non-negative values on the other vertices.
Thus, $h_x$ naturally induces a well-defined function on ${\rm ind}({\mathcal{C}}_Q)$,
which we still denote as $h_x$ on ${\rm ind}({\mathcal{C}}_Q)$. We extend $h_x$ to a function defined on ${\mathcal{C}}_Q$
by requiring that $h_x(X \oplus Y) = h_x(X) + h_x(Y)$ for all objects $X, Y$ of ${\mathcal{C}}_Q$.
Let ${\mathcal{S}}'_x$ be the slice in ${\mathbb{Z}}Q$ with $x$ its unique sink and ${\mathcal{S}}''_x$ the slice in ${\mathbb{Z}}Q$ with $x$ its unique source.

Let $Z$ be an indecomposable object of ${\mathcal{C}}_Q$. If there is an arrow from $x$ to $Z$ in the Auslander-Reiten quiver of ${\mathcal{C}}_Q$,
then $Z$ and $\tau Z$ both lie in the $(\tau^{-1} \Sigma)$-orbit of the convex hull of ${\mathcal{S}}'_x$ and ${\mathcal{S}}''_x$.
Thus, both $h_x(Z)$ and $h_x(\tau Z)$ are zero,
which impies that all $h_x(y)$ appearing in the right hand side of item $c)$ are non-positive. Hence, we have that
$$h_x(Z) + h_x(\tau Z) = \sum_{y \rightarrow Z} \max \{h_x(y), 0 \} = 0 =  \max \{\sum_{y \rightarrow Z} h_x(y), 0 \},$$
where `$y \rightarrow Z$' are arrows in $\Gamma$.
If there is no arrow from $x$ to $Z$ in the Auslander-Reiten quiver of ${\mathcal{C}}_Q$, then we have the following equalities
\begin{align*}
h_x(Z) + h_x(\tau Z) & = \sum_{y \rightarrow Z} \max \{h_x(y), 0 \} \\
& = \sum_{y \rightarrow Z} h_x(y) =  \max \{\sum_{y \rightarrow Z} h_x(y), 0 \},
\end{align*}
where `$y \rightarrow Z$' are arrows in $\Gamma$.
%from item $c)$ above we can read that all $h_x(y)$ are zero where $y$ is an indecomposable direct summand of the middle $E$ appearing in the non-split triangle
%$(\ast)$ before Proposition \ref{4} associated to $Z$
Therefore, for all non-split triangles as the triangles $(\ast)$ before Proposition \ref{4},
the function $h_x$ satisfies item d3) in Definition \ref{1}. Besides, by Proposition \ref{16}, there is a unique tropical frieze
$g: {\mathcal{C}}_Q \rightarrow {\mathbb{Z}}$ such that $g(Y) = h_x(Y)$
for all indecomposables $Y$ which come from the same slice containing $x$.
Thus, we have that $h_x = g$ and $h_x$ is a tropical frieze on ${\mathcal{C}}_Q$.

Let ${\mathcal{S}}_x$ be any slice in $\mathbb{Z}$Q with $x$ a source. Set $T = \oplus_{Y \in {\mathcal{S}}_x} Y$. It is a basic cluster-tilting object of ${\mathcal{C}}_Q$.
Let $B$ be the endomorphism algebra of $T$ and
$S_x$ the simple $B$-module corresponding to $x$. Clearly ${\mathcal{S}}_x$ is the quiver of $B$. Set $m = -S_x$.
Then $f_{T, -S_x}$ is a tropical frieze and takes the same values as $h_x$ on all indecomposable direct summands of $T$. As a result, the function
$h_x$ is equal to $f_{T, -S_x}$.

However, the sum $\sum_{x}h_x$ of cluster-hammock functions with
all $x$ coming from the same slice $\mathcal{S}$ in ${\mathbb{Z}}Q$ is not always a tropical frieze,
which is quite different to the Corollary in Section 6 of \cite{Rin11}. Here we also use the same counter-example on ${\mathcal{C}}_{\vec{A_3}}$
as in Section 2. We already know that the functions $d_1$ and $h_{\Sigma P_1}$ are tropical friezes. Let $T = \Sigma P_1 \oplus \Sigma P_2 \oplus \Sigma P_3$.
Then $d_1$ and $h_{\Sigma P_1}$ coincide on all $\Sigma P_i \, (1 \leq i \leq 3)$. Thus, $h_{\Sigma P_1}$ is equal to $d_1$. Similarly, the tropical frieze
$h_{\Sigma P_3}$ is equal to $d_3$. But the sum $h_{\Sigma P_1} + h_{\Sigma P_3} = d_1 + d_3$ is not a tropical frieze.

\vspace{.3cm}

\section{simple illustrations for the cases $A_1$ and $A_2$}

Let us first look at the cluster category ${\mathcal{C}} = {\mathcal{C}}_Q$ of the quiver $Q$ of type $A_1$. Let $X$ and $\Sigma X$ be
the two indecomposable objects in ${\mathcal{C}}_{Q}$. Assume $f$ is a tropical frieze on ${\mathcal{C}}_{Q}$. Then we have that
$$f(X) + f(\Sigma X) = 0.$$ Set $T= X$ and $m = f(X)S_X$, where $S_X$ is the unique simple $({\rm End}_{{\mathcal{C}}_{Q}}(X))$-module.
Since $\langle S_X, m \rangle_a$ is zero, by Theorem \ref{2} the function $f_{T,m}$ is a tropical frieze.
The following equalities
\begin{align*}
f_{T,m}(X) & = \langle F({\rm ind}(X)),f(X) S_X \rangle = f(X) \quad \mbox{and} \\
f_{T,m}(\Sigma X) & = \langle F({\rm ind}(\Sigma X)),f(X) S_X \rangle = - f(X) = f(\Sigma X)
\end{align*}
clearly hold. Therefore, the tropical frieze $f$ is equal to $f_{T,m}$.

Now let us look at the cluster category ${\mathcal{C}} = {\mathcal{C}}_Q$ of a quiver $Q$ of type $A_2$. Assume that $f$ is a non-zero tropical frieze on
${\mathcal{C}}_Q$. Following Proposition \ref{4}, we know that there exist an indecomposable object $X$ such that $f(X) < 0$.
Let $Y$ and $Y'$ be the two non-isomorphic indecomposables such that $X \oplus Y$ and $X \oplus Y'$ are cluster-tilting objects of ${\mathcal{C}}_Q$. Then we have that
$$f(Y) + f(Y') = \max \{f(X), 0 \} = 0. $$ Therefore, there must
 exist a cluster-tilting object $T = T_1 \oplus T_2$ with $T_i$ indecomposable
such that $$f(T_1) \geq 0 \quad \mbox{and} \quad f(T_2) < 0.$$ Let $Q_T$ be the quiver of the endomorphism algebra
$B = {\rm End}_{{\mathcal{C}}_Q}(T)$. Let $P_i$ be the indecomposable projective $B$-module and $S_i$ its corresponding simple top. The quiver $Q_T$ is also
of type $A_2$.

If $S_1$ attaches to the sink in $Q_T$, set $m = f(T_1) S_1 + f(T_2)S_2$, then
\begin{align*}
\langle S_1, m \rangle_a & = -f(T_2) {\rm dim} {\rm Ext}^1_B (S_1, S_2) = -f(T_2) > 0 \quad \mbox{and} \\
\langle S_2, m \rangle_a & = f(T_1) {\rm dim} {\rm Ext}^1_B (S_1, S_2) = f(T_1) \geq 0,
\end{align*}
which implies that $f_{T,m}$ is a tropical frieze by Theorem \ref{2}. Moreover, the tropical friezes $f$ and $f_{T,m}$ coincide on $T_i$. Therefore, the
tropical frieze $f$ is equal to $f_{T,m}$.

If $S_1$ attaches to the source in $Q_T$, set $T' = \mu_2 \mu_1 (T) = T'_1 \oplus T'_2$,
where $T'_1$ and $T'_2$ come from the following non-split triangles in ${\mathcal{C}}_Q$
$$T_1 \rightarrow T_2 \rightarrow T'_1 \rightarrow \Sigma T_1, \qquad T'_1 \rightarrow 0 \rightarrow T_1 \rightarrow \Sigma T'_1;$$
$$T_2 \rightarrow T'_1 \rightarrow T'_2 \rightarrow \Sigma T_2, \qquad T'_2 \rightarrow 0 \rightarrow T_2 \rightarrow \Sigma T'_2.$$
We can calculate that
$$f(T'_1) = -f(T_1) \leq 0 \quad \mbox{and} \quad f(T'_2) = -f(T_2) > 0.$$ Notice that the quiver $Q_{T'}$ of the endomorphism algebra
$B' = {\rm End}_{{\mathcal{C}}_Q}(T')$ is $T'_1 \rightarrow T'_2$. Let $S'_i$ be the simple $B'$-module corresponding to $T'_i$.
Now we go back to the above cases.
Set $m' = f(T'_1) S'_1 + f(T'_2)S'_2$.
Then we have that $f_{T',m'}$ is a tropical frieze and takes the same values as $f$ on $T'_i$. Thus, the tropical frieze $f$ is equal to $f_{T',m'}$.

In fact, such a phenomenon for the cases $A_1$ and $A_2$ is a common phenomenon for the Dynkin case, which we will state in Theorem \ref{15} in the next section.

Let $f_{T,m}$ be a tropical frieze on ${\mathcal{C}}_Q$ with $Q$ a quiver of type $A_2$. Suppose that the quiver $Q_T$ of the
endomorphism algebra $B = {\rm End}_{{\mathcal{C}}_Q}(T)$ is $(T_1 \rightarrow T_2)$ and $m = m_1S_1 +m_2S_2$. From Remark \ref{3}
we know that $\langle S_i, m \rangle_a \geq 0$ for $i = 1, 2$, that is,
\begin{align*}
\langle S_1, m \rangle_a & = m_2 {\rm dim} {\rm Ext}^1_B (S_2, S_1) = m_2 \geq 0, \quad \mbox{and} \\
\langle S_2, m \rangle_a & = -m_1 {\rm dim} {\rm Ext}^1_B (S_2, S_1) = -m_1 \geq 0.
\end{align*}
Notice that $f_{T,m}(T_i) = \langle FT_i, m \rangle = m_i$ for $i
=1, 2$. Set $T' = \mu_1(T) = T'_1 \oplus T_2$ and $T'' = \mu_2(T) =
T_1 \oplus T''_2$. Then the following expressions hold
\begin{align*}
f_{T,m}(T'_1) & = \max \{ f_{T,m}(T_2), 0 \} - f_{T,m}(T_1) \geq -f_{T,m}(T_1) \geq 0, \\
f_{T,m}(T''_2) & = \max \{ f_{T,m}(T_1), 0 \} - f_{T,m}(T_2) =
-f_{T,m}(T_2) \leq 0.
\end{align*}
Therefore, in the $A_2$ case, there exist cluster-tilting objects $T'$ and $T''$ such that
$f_{T,m}$ takes non-negative values on direct summands of $T'$ and non-positive values on direct summands of $T''$.

\vspace{.3cm}

\section{the main theorem (dynkin case)}

As a generalization of the phenomenon illustrated in Section 4, the aim of this section is to show the following theorem:

\begin{thm}\label{15}
Let ${\mathcal{C}}_Q$ be the cluster category of a Dynkin quiver $Q$. Then all tropical friezes on ${\mathcal{C}}_Q$ are of the form $f_{T,m}$,
where $T$ is a cluster-tilting object and $m$ an element in the Grothendieck group $K_0({\rm mod} {\rm End}_{{\mathcal{C}}_Q}(T))$.
\end{thm}

We will prove the theorem in sections 5.1 and 5.2. First, we need to
introduce some notation:

Let $\mathcal{C}$ be a 2-Calabi-Yau category with cluster-tilting
object. Let $f$ be a tropical frieze on the category $\mathcal{C}$
and $T = T_1 \oplus \ldots \oplus T_n$ a basic cluster-tilting
object of $\mathcal{C}$. Suppose that the quiver $Q$ of the
endomorphism algebra ${\rm End}_{\mathcal{C}}(T)$ does not have
loops nor $2$-cycles. Let $b_{ij}$ denote the number of arrows $i
\rightarrow j$ minus the number of arrows $j \rightarrow i$ in $Q$
(notice that at least one of these two numbers is zero). For each
integer $1 \leq i \leq n$, let $g_i(T)$ be the integer
\begin{displaymath}
g_i(T) =\sum_{r}[b_{ri}]_+ f(T_r) - \sum_s [b_{is}]_+ f(T_s),
\end{displaymath}
where $[b_{kl}]_+ = \max\{b_{kl}, 0\}$ is equal to the number of arrows $k \rightarrow l$ in $Q$.
Denote by $g(T)$ the class
$\sum^n_{i=1} g_i(T)[T_i]$ in the Grothendieck group $K_0 ({\rm add}T)$.

\subsection{Transformations of the class $g(T)$ under mutations}

Since the quiver $Q$ does not have loops, for each $T_k$, there is a unique indecomposable object $T'_k$ such that
the space ${\rm Ext}^1_{\mathcal{C}}(T'_k, T_k)$ is one-dimensional and the non split triangles are given \cite{Ke10} by
$$T'_k \rightarrow E \rightarrow T_k \rightarrow \Sigma T'_k \quad \mbox{and} \quad T_k \rightarrow E' \rightarrow T'_k\rightarrow \Sigma T_k,$$
where $$E =\bigoplus_r [b_{rk}]_+ T_r \quad \mbox{and} \quad  E'=\bigoplus_s [b_{ks}]_+ T_s.$$
Let $T' = \mu_k (T) = T'_k \oplus (\bigoplus_{i \neq k} T_i)$.
Define linear transformations $\phi_+$ and $\phi_{-}$ from $K_0({\rm add}T)$ to $K_0({\rm add}T')$ as in \cite{DK} by
\begin{align*}
\phi_+ (T_i) & = \phi_{-} (T_i) = [T_i] \quad \quad \mbox{for}\, i \neq k, \, \mbox{and}\\
\phi_+ (T_k) & = [E]-[T'_k] = -[T'_k] + \sum_r [b_{rk}]_+[T_r] \\
\phi_{-} (T_k) & = [E'] - [T'_k] = -[T'_k] + \sum_s [b_{ks}]_+ [T_s].
\end{align*}
It was shown in \cite{DK} that if $X$ is a rigid object of $\mathcal{C}$, then the index of $X$ with respect to cluster-tilting objects transforms as follows:
\begin{displaymath}
{\rm ind}_{T'}(X) = \left\{ \begin{array}{ll}
       \phi_+({\rm ind}_{T}(X)) & \textrm{if}\,\, [{\rm ind}_T(X) : T_k] \geq 0,\\
       \phi_{-}({\rm ind}_{T}(X)) & \textrm{if}\, \, [{\rm ind}_T(X) : T_k] \leq 0,
\end{array} \right.
\end{displaymath}
where $[{\rm ind}_T(X) : T_k]$ denotes the coefficient of $T_k$ in the decomposition of ${\rm ind}_T (X)$ in the category $K_0(\mbox{add}T)$.

\begin{prop} \label{10}
Suppose that the quivers $Q$ and $Q'$ of the endomorphism algebras ${\rm End}_{\mathcal{C}}(T)$ and ${\rm End}_{\mathcal{C}}(T')$ do not
have loops nor 2-cycles.
Then the element $g(T)$ transforms in the same way as above, i.e.
\begin{displaymath}
g(T') = \left\{ \begin{array}{ll}
       \phi_+(g(T)) & \textrm{if}\,\, g_k(T) \geq 0,\\
       \phi_{-}(g(T)) & \textrm{if}\,\, g_k(T) \leq 0.
\end{array} \right.
\end{displaymath}
\end{prop}
\begin{pf}
We first assume that $g_k(T) \geq 0$, that is,
$f(E) =\sum_r [b_{rk}]_+ f(T_r) \geq \sum_s [b_{ks}]_+ f(T_s) = f(E').$
Since $f$ is a tropical frieze, we have that $f(T_k) + f(T'_k) = f(E) =\sum_r [b_{rk}]_+ f(T_r).$ We compute $\phi_+(g(T))$:
\begin{align*}
\phi_+(g(T)) & = \phi_+(\sum^n_{i=1} g_i(T)[T_i]) = \sum_{i \neq k} g_i(T)[T_i] + g_k(T) \phi_+(T_k) \\
& =\sum_{i \neq k} g_i(T)[T_i] - g_k(T)[T'_k] + \sum_{r} g_k(T) [b_{rk}]_+ [T_r] \\
& = \sum_{i \neq k} (g_i(T)+[b_{ik}]_+g_k(T))[T_i] - g_k(T)[T'_k].
\end{align*}

By assumption, the quivers $Q$ and $Q'$ do not have loops nor 2-cycles. Following \cite{BIRSc},
we know that $Q' = \mu_k (Q)$ is the mutation of the quiver $Q$ at vertex $k$. Let $b'_{ij}$ denote the number of arrows
$i \rightarrow j$ minus the number of arrows $j \rightarrow i$ in $Q'$. Then it is known from \cite{FZI} that
\begin{displaymath}
b'_{ij} = \left\{ \begin{array}{ll}
       b_{ji} & \textrm{if}\, i = k \, \textrm{or}\, j = k,\\
       b_{ij}+\frac{|b_{ik}|b_{kj}+b_{ik}|b_{kj}|}{2} & \textrm{otherwise}.
\end{array} \right.
\end{displaymath}

It is obvious that
\begin{align*}
g_k(T') & =  \sum_{r} [b'_{rk}]_+ f(T_r)\, - \sum_{s} [b'_{ks}]_+ f(T_s) \\
& = \sum_{r} [b_{kr}]_+ f(T_r)\, - \sum_{s} [b_{sk}]_+ f(T_s) = -g_k(T).
\end{align*}

For vertices $i \neq k$, we distinguish three cases to compute $g_i(T')$.

If $b_{ik} = b_{ki} = 0$, then $b_{ij}' = b_{ij}$ and $b'_{ji} = b_{ji}$ for all vertices $j$. In this case, we have that
\begin{align*}
 g_i (T') & = \sum_{r} [b'_{ri}]_+ f(T_r) - \sum_{s} [b'_{is}]_+f(T_s) \\
& = \sum_{r} [b_{ri}]_+ f(T_r) -  \sum_{s} [b_{is}]_+ f(T_s) = g_i (T).
\end{align*}

If $b_{ik} > 0$, then
\begin{align*}
g_i (T') & = \sum_{r} [b'_{ri}]_+ f(T_r) - \sum_{s} [b'_{is}]_+ f(T_s) = (\sum_{r} [b_{ri}]_+ f(T_r) + b_{ik}f(T'_k)) \\
 & \qquad \qquad - (\sum_{s} [b_{is}]_+ f(T_s) - b_{ik}f(T_k) + \sum_{s'} b_{ik} [b_{ks'}]_+ f(T_{s'}))\\
& = g_i(T) + b_{ik} (f(T'_k)+ f(T_k) -\sum_{s'} [b_{ks'}]_+ f(T_{s'})) \\
& = g_i(T) + b_{ik} (\sum_{r} [b_{rk}]_+ f(T_r) -\sum_{s} [b_{ks}]_+ f(T_{s})) \\
& = g_i(T) + b_{ik}g_k(T).
\end{align*}

If $b_{ik} < 0$, then $b_{ki} = -b_{ik} > 0$, and
\begin{align*}
g_i(T') & = \sum_{r} [b'_{ri}]_+ f(T_r) - \sum_{s} [b'_{is}]_+ f(T_s) \\
& =(\sum_{r} [b_{ri}]_+ f(T_r) - b_{ki}f(T_k) + \sum_{r'} [b_{r'k}]_+ b_{ki}f(T_{r'})) \\
& \qquad \qquad \qquad - (\sum_s [b_{is}]_+ f(T_s) + b_{ki}f(T'_k)) \\
& = g_i(T) - b_{ki} (f(T_k)+ f(T'_k) -\sum_{r'} [b_{r'k}]_+ f(T_{r'})) \\
& = g_i(T) - b_{ki} (\sum_{r} [b_{rk}]_+ f(T_r) -\sum_{r'} [b_{r'k}]_+ f(T_{r'})) = g_i(T).
\end{align*}

Therefore, we obtain that $g(T') = \phi_+(g(T))$ when $g_k(T) \geq 0$. In a similar way we can also obtain that $g(T') = \phi_{-}(g(T))$
when $g_k(T) \leq 0$.
\end{pf}

\subsection{Proof of the main theorem}

Let $T_0$ and $T_1$ be two objects in add$T$ which do not have a direct summand in common. Let $\eta$ be a morphism in ${\mathcal{C}}(T_1,T_0)$.
Denote by $C(\eta)$ the cone of $\eta$. Then we have the following triangle in $\mathcal{C}$
$$\qquad \qquad T_1 \stackrel{\eta}\rightarrow T_0 \rightarrow C(\eta) \rightarrow \Sigma T_1. \qquad \qquad \qquad \qquad (\ast \ast)$$
The algebraic group $\mbox{Aut}(T_0) \times \mbox{Aut}(T_1)$ acts on ${\mathcal{C}}(T_1,T_0)$ via
$$(g_0,g_1){\eta}' = g_0 {\eta}' g_1^{-1}.$$
Let ${\mathcal{O}}_{\eta}$ denote the orbit of $\eta$ in the space ${\mathcal{X}} := {\mathcal{C}}(T_1,T_0)$
under the above action of $\mbox{Aut}(T_0) \times \mbox{Aut}(T_1)$.

It is not hard to obtain the following lemma. For the convenience of the reader we include a proof.
\begin{lem} \label{7}
Let $\eta$ and $\eta'$ be two morphisms in $\mathcal{X}$. Then ${\mathcal{O}}_{\eta} = {\mathcal{O}}_{\eta'}$
if and only if $C(\eta) \simeq C(\eta')$.
\end{lem}

\begin{pf}
First we assume that ${\mathcal{O}}_{\eta} = {\mathcal{O}}_{\eta'}$.
Then there exists an element $(g_0,g_1) \in \mbox{Aut}(T_0) \times \mbox{Aut}(T_1)$ such that
$\eta' = g_0 \eta g_1^{-1}$. The commutative square $g_0 \eta = \eta' g_1$ can be completed to a commutative diagram of triangles as follows
\[
\xymatrix{ T_1 \ar[r]^{\eta}
\ar[d]^{g_1} & T_0 \ar[r]^{\iota} \ar[d]^{g_0} & C(\eta) \ar[r]^{p} \ar[d]^h & \Sigma T_1
\ar[d]^{\Sigma g_1} \\
T_1 \ar[r]^{\eta'} & T_0 \ar[r]^{\iota'}& C(\eta') \ar[r]^{p'} & \Sigma T_1.}
\]
Here the morphism $h$ is an isomorphism from $C(\eta)$ to $C(\eta')$.

Second we assume that $C(\eta) \simeq C(\eta')$. Let $h$ be an isomorphism from $C(\eta)$ to $C(\eta')$ and $h^{-1}$ its inverse.
Since the space ${\mathcal{C}}(T_0, \Sigma T_1)$ vanishes, we have that (keeping the notation as in the above commutative diagram)
$$p' h \iota = 0 \quad \mbox{and} \quad p h^{-1} \iota' = 0. $$
Thus, there exist two morphisms $g_0$ and $g'_0$ in ${\mathcal{C}}(T_0,T_0)$ such that
$$\iota' g_0 = h \iota \quad \mbox{and} \quad \iota g'_0 = h^{-1} \iota'.$$
As a consequence, the equalities
$$\iota g'_0 g_0 = h^{-1} \iota' g_0 = h^{-1} h \iota = \iota \quad \mbox{and} \quad \iota' g_0 g'_0 = h \iota g'_0 = hh^{-1} \iota' = \iota'$$
hold. Thus, we have that $g_0 g'_0 = 1 = g'_0 g_0$. The morphism $g_0$ is an element in $\mbox{Aut}(T_0)$. The commutative square $h \iota = \iota' g_0$
can be completed to a commutative diagram of triangles as above. Thhus, there exists an element $g_1 \in \mbox{Aut}(T_1)$ such that $g_0 \eta = \eta' g_1$.
Therefore, the two orbits ${\mathcal{O}}_{\eta}$ and ${\mathcal{O}}_{\eta'}$ are the same.
\end{pf}

\begin{lem}\label{8}
Keep the above notation. We have the equality
$${\rm codim}_{\mathcal{X}} {\mathcal{O}}_{\eta} = 1/2\, {\rm dim} {\rm Ext}^1_{\mathcal{C}}(C(\eta),C(\eta)).$$
\end{lem}

\begin{pf}
Let $F$ be the functor ${\mathcal{C}}(T,?)$ and $B$ the endomorphism algebra of $FT$. We denote the space ${\rm Hom}_B (FT_1, FT_0)$ by $F{\mathcal{X}}$.
Since $F$ induces a category equivalence from ${\mathcal{C}}/{\rm add}(\Sigma T)$ to ${\rm mod}B$, we have that
$$\mbox{codim}_{\mathcal{X}} {\mathcal{O}}_{\eta} = \mbox{codim}_{F{\mathcal{X}}} {\mathcal{O}}_{F{\eta}}.$$

The algebra $B$ is a finite-dimensional algebra, both $FT_1$ and $FT_0$ are finitely generated $B$-modules. As in \cite{Pla11}, we
view $F{\eta}$ as a complex in $K^b(\mbox{proj} B)$ and
define the space $E(F{\eta})$ as
$$E(F{\eta}) = {\rm Hom}_{K^b (\mbox{proj} B)} (\Sigma^{-1} F{\eta}, F{\eta}).$$
Following Lemma 2.16 in \cite{Pla11}, we have the equality $$\mbox{codim}_{F{\mathcal{X}}} {\mathcal{O}}_{F{\eta}} = {\rm dim} E(F{\eta}).$$
The exact sequence $$FT_1 \stackrel{F{\eta}}\rightarrow FT_0 \rightarrow F(C({\eta})) \rightarrow 0$$ is a minimal projective presentation of $F(C({\eta}))$.
Still following from \cite{Pla11}, the equality
$${\rm dim} E(F{\eta}) = {\rm dim} {\rm Hom}_B (F(C({\eta})), \tau F(C({\eta})))$$
holds, where $\tau$ is the Auslander-Reiten translation. Moreover, by Section 3.5 in \cite{KR07}, we have that $F(\Sigma C({\eta})) \simeq \tau F(C({\eta}))$.

For two objects $X$ and $Y$ of $\mathcal{C}$, let $(\Sigma T)(X,Y)$ be the subspace of ${\mathcal{C}}(X,Y)$ consisting of morphisms from $X$ to $Y$
factoring through an object in add$(\Sigma T)$, let ${\mathcal{C}}/_{(\Sigma T)} (X, Y)$ denote the space ${\mathcal{C}}(X,Y)/{(\Sigma T)(X,Y)}$.
Lemma 3.3 in \cite{Pal08} shows that there is a bifunctorial isomorphism
$${\mathcal{C}}/_{(\Sigma T)} (X, \Sigma Y) \simeq D (\Sigma T) (Y, \Sigma X).$$
If we choose $Y = X$, then we can deduce that
$${\rm dim}\, {\mathcal{C}}/_{(\Sigma T)} (X, \Sigma X)  = {\rm dim} (\Sigma T) (X, \Sigma X)  = 1/2 \, {\rm dim}\,{\mathcal{C}} (X, \Sigma X).$$
Notice that the equivalence $F$ gives the following equality
$${\rm dim} {\rm Hom}_B (F(C({\eta})), \tau F(C({\eta}))) = {\rm dim}\, {\mathcal{C}}/_{(\Sigma T)} (C({\eta}), \Sigma C({\eta})).$$
Finally, if we combine all the equalities about dimensions together, then we can obtain that
$$\mbox{codim}_{\mathcal{X}} {\mathcal{O}}_{\eta} = 1/2\, {\rm dim} {\rm Ext}^1_{\mathcal{C}}(C({\eta}), C({\eta})).$$
\end{pf}

If we do not assume that $T_0$ and $T_1$ do not have a common direct summand, then the equality in Lemma \ref{8} becomes
$$\mbox{codim}_{\mathcal{X}} {\mathcal{O}}_{\eta} \geq 1/2\, {\rm dim} {\rm Ext}^1_{\mathcal{C}}(C({\eta}),C({\eta})).$$
This is because the third equality in the proof becomes
$${\rm dim} E(F{\eta}) \geq {\rm dim} {\rm Hom}_B (F(C({\eta})), \tau F(C({\eta})))$$ for arbitrary projective presentations.

\begin{lem} \label{6}
Suppose that $\mathcal{C}$ has only finitely many isomorphism classes of indecomposable objects. Then the set
$\{  [C({\eta})] | {\eta} \in {\mathcal{C}}(T_1,T_0)
\}$ is finite, where $[C({\eta})]$ denotes the isomorphism class of $C({\eta})$ in $\mathcal{C}$.
\end{lem}

\begin{pf}
We use the same exact sequence
$$FT_1 \stackrel{F{\eta}}\rightarrow FT_0 \rightarrow F(C({\eta})) \rightarrow 0$$
as in the proof of Lemma \ref{8},
which is a projective presentation of $F(C({\eta}))$. By assumption $\mathcal{C}$ has only finitely many isomorphism classes of indecomposable objects.
So the number of isomorphism classes of indecomposable $B$-modules is also finite.
Notice that the dimension of $F(C({\eta}))$ is bounded by the dimension of $FT_0$.
Hence, the set of $\{  [F(C({\eta}))] | {\eta} \in {\mathcal{C}}(T_1,T_0)
\}$ is finite, where $[F(C({\eta}))]$ denotes the isomorphism class of $F(C({\eta}))$ in ${\rm mod}B$.

Now we decompose $C({\eta})$ as $X_{\eta} \oplus \Sigma T_{\eta}$, where $X_{\eta}$ does not contain a direct summand in add$T$. We have that $F(C({\eta})) = F(X_{\eta})$.
Since ${\mathcal{C}}(T_0, \Sigma T_{\eta})$ vanishes,
we can rewrite the triangle $(\ast \ast)$ before Lemma \ref{7} as
$$T_1 \stackrel{
\left(\begin{array}{c}
{\eta} \\
0
\end{array}\right)
} \longrightarrow T_0 \oplus 0 \stackrel{
\left(\begin{array}{cc}
\iota_{\eta} & 0 \\
0 & 0
\end{array}\right)
} \longrightarrow X_{\eta} \oplus \Sigma T_{\eta} \rightarrow \Sigma T_1,$$
which is the direct sum of the following two triangles
$$\Sigma^{-1} C(\iota_{\eta}) \rightarrow T_0 \stackrel{\iota_{\eta}}\rightarrow X_{\eta} \rightarrow C(\iota_{\eta}), \qquad \mbox{and}$$
$$T_{\eta} \rightarrow 0 \rightarrow \Sigma T_{\eta} \rightarrow \Sigma T_{\eta}. \qquad \quad$$ Here $C(\iota_{\eta})$ denotes the cone of the morphism $\iota_{\eta}$.
Therefore, the object $T_{\eta}$ is a direct summand of $T_1$, and there are only finitely many choices. In conclusion, there are only finitely many
isomorphism classes of $C({\eta})$ when ${\eta}$ runs over the space ${\mathcal{C}}(T_1,T_0)$.
\end{pf}

Under the assumption that $\mathcal{C}$ has only finitely many isomorphism classes of indecomposable objects,
by combining Lemma \ref{7} and Lemma \ref{6} we can obtain that there are only finitely many orbits ${\mathcal{O}}_{\eta}$ in the affine space $\mathcal{X}$.
Therefore, there must exist some morphism ${\eta}$ such that
$$\mbox{codim}_{\mathcal{X}} {\mathcal{O}}_{\eta} = 0,$$ which implies that $C({\eta})$ is a rigid object by Lemma \ref{8}.
We say a morphism $\eta$ {\it generic} if its cone $C(\eta)$ is rigid.
We deduce the following proposition

\begin{prop}\label{9}
Suppose that $\mathcal{C}$ has only finitely many isomorphism classes of indecomposable objects.
Then there exists a generic morphism
${\eta} \in {\mathcal{C}}(T_1,T_0)$ with the cone $C({\eta})$ rigid.
\end{prop}

Now we are ready to prove our main theorem.

\smallskip

{\em{Proof of Theorem \ref{15}.}}
Let $T = T_1 \oplus \ldots \oplus T_n$ be any basic cluster-tilting object in ${\mathcal{C}}_Q$. Keeping the notation in the beginning of this section,
we define two objects
$$ L =\bigoplus_{ g_i(T)<0 } T_i^{ -g_i (T) }  \qquad \mbox{and} \qquad R =\bigoplus_{ g_i(T)>0 } T_i^{ g_i (T) }. $$
By Proposition \ref{9}, there exists a morphism ${\eta} \in {\rm Hom}_{{\mathcal{C}}_Q}(L,R)$ such that the cone $C({\eta})$ is rigid. The triangle
$$ L \stackrel{{\eta}}\rightarrow R \rightarrow C({\eta}) \rightarrow \Sigma L$$ implies that the index
$$\mbox{ind}_T (C({\eta})) = [R] - [L] = g(T).$$

Since $C({\eta})$ is rigid, there exists a cluster-tilting object $T'$ of ${\mathcal{C}}_Q$ such that $C({\eta}) \in {\rm add}T'$. The triangle
$$\Sigma^{-1} C({\eta}) \rightarrow 0 \rightarrow C({\eta}) \rightarrow C({\eta})$$ gives us that
$${\rm ind}_{\Sigma^{-1} T'} (C({\eta})) \in {\mathbb{Z}}^n_{\leq 0}.$$
Set $T'' = \Sigma^{-1} T'$. It was shown in \cite{BMR08} that
the quiver of the endomorphism algebra of a cluster-tilting object of ${\mathcal{C}}_Q$
does not have loops nor 2-cycles. Therefore, it
follows from Proposition \ref{10} that
$$g(T'') = {\rm ind}_{T''} (C({\eta})) \in {\mathbb{Z}}^n_{\leq 0},$$
that is, $g_i(T'') \leq 0$.

Let $B''$ denote the endomorphism algebra ${\rm End}_{{\mathcal{C}}_Q} (T'')$ and $Q''$ its associated quiver.
Let $S^{''}_i$ be the simple top of the indecomposable projective $B''$-module
$P''_i = {\rm Hom}_{{\mathcal{C}}_Q}(T'', T''_i)$. Set $$m'' = \sum^n_{j=1} f(T''_j) S''_j \, (\in K_0({\rm mod}B'')).$$
Then for each simple $B''$-module $S''_i$,
we have that
\begin{align*}
\langle S''_i , m'' \rangle_a & = \sum^n_{j=1} f(T''_j) \langle S''_i , S''_j \rangle_a \\
& =  \sum^n_{j=1} f(T''_j) (-{\rm dim} {\rm Ext}^1_{B''} (S''_i, S''_j)) +
\sum^n_{j=1} f(T''_j) {\rm dim} {\rm Ext}^1_{B''} (S''_j, S''_i) \\
& = -\sum^n_{j=1} [b''_{ji}]_+ f(T''_j) + \sum^n_{j=1} [b''_{ij}]_+ f(T''_j) = -g_i(T'') \geq 0,
\end{align*}
where $b''_{kl}$ denotes the number of arrows
$k \rightarrow l$ minus the number of arrows $l \rightarrow k$ in $Q''$.
Therefore, by Theorem \ref{2} the function $f_{T'',m''}$ is a tropical frieze.
Since we have
$$f_{T'',m''}(T''_i) = \langle P''_i, m'' \rangle = \langle P''_i, f(T''_i) S''_i  \rangle = f(T''_i),$$
the tropical friezes $f$ and $f_{T'',m''}$ coincide on all $T''_i$. Now it follows from Proposition \ref{11} that $f$ is equal to $f_{T'',m''}$.

\subsection{Sign-coherence property}

For any tropical frieze $f$ on ${\mathcal{C}}_Q$ with $Q$ a Dynkin quiver,
we will see in this subsection the existence of cluster-tilting objects whose indecomposable direct summands have sign-coherent values under $f$.

\begin{thm}\label{17}
Let ${\mathcal{C}}_Q$ be the cluster category of a Dynkin quiver $Q$ and $f$ a tropical frieze on ${\mathcal{C}}_Q$. Then there
exists a cluster-tilting object $T$ such that $$f(T_i) \geq 0  \quad (\mbox{resp.}\,\, f(T_i) \leq 0 )$$ for
all indecomposable direct summands $T_i$ of $T$.
%and also exists a cluster-tilting object $T''$ such that $f(T''_i) \geq 0$ for all indecomposable direct summand $T''_i$ of $T''$.
\end{thm}

\begin{pf}
Since $f$ is a tropical frieze on ${\mathcal{C}}_Q$, it follows from Theorem \ref{15} that $f$ is equal to some $f_{T,m}$
with $T$ a cluster-tilting object and $m$ an element in $K_0 ({\rm mod} {\rm End}_{{\mathcal{C}}_Q}(T))$. We divide the proof
into three steps.

Step $1$. For any cluster-tilting object $S$ of ${\mathcal{C}}_Q$, we define its associated positive cone as
$$C(S) = \{ \mbox{ind}_T(U) |\, U \in {\rm add}S \} \, (\subset K_0 ({\rm add}T)).$$ Each element $X \in K_0 ({\rm add}T)$ can be
written uniquely as
$$X = [T_0] - [T_1],$$ where $T_0$, $T_1 \in {\rm add}T$ without common indecomposable direct summands. By Proposition \ref{9},
there exists some morphism $\eta \in {\rm Hom}_{{\mathcal{C}}_Q}(T_1,T_0)$ such that the cone $C(\eta)$ is rigid. Moreover, we
have that $${\rm ind}_T (C(\eta)) = [T_0] - [T_1] = X.$$ Since $C(\eta)$ is rigid, it belongs to add$S$ for some cluster-tilting
object $S$ of ${\mathcal{C}}_Q$, which implies that the element $X$ belongs to the positive cone $C(S)$. As a consequence, we can
obtain that $$K_0 ({\rm add}T) = \bigcup_S C(S),$$ where $S$ ranges over all (finitely many) cluster-tilting objects of ${\mathcal{C}}_Q$.

Step $2$. Let $T_1, \ldots, T_n$ be the pairwise non-isomorphic indecomposable direct summands of $T$. Suppose that
$m=\sum^n_{i=1} m_i S_i$ with $S_i$ the simple ${\rm End}_{{\mathcal{C}}_Q}(T)$-module corresponding to $T_i$. Set
$$H_m^{\geq 0} = \{ X \in K_0({\rm add}T) \,  | \, \langle FX, m \rangle \geq 0 \}.$$ It is clear that
$$\langle {\rm sgn}(m_i) FT_i, m \rangle = |m_i| \geq 0,$$ where $F$ is the functor ${\rm Hom}_{{\mathcal{C}}_Q}(T, ?)$ and
\begin{displaymath}
{\rm sgn}(m_i) = \left\{ \begin{array}{ll}
       1 & \textrm{if}\, m_i \geq 0,\\
       -1 & \textrm{if}\, m_i < 0.
\end{array} \right.
\end{displaymath}
Let $H$ be the hyperquadrant of $K_0({\rm add}T)$ consisting of the non-negative linear combinations of the ${\rm sgn}(m_i)[T_i],\,1 \leq i \leq n$.
Then we have that $$H \subset H_m^{\geq 0}.$$

Step $3$. It was shown in Section 2.4 of \cite{DK} that each positive cone $C(S)$
is contained in a hyperquadrant of $K_0({\rm add}T)$ with respect to the given basis $[T_i],\, 1 \leq i \leq n$.
Thus, each hyperquadrant of $K_0({\rm add}T)$ is a union of positive cones. Let $T'$ be a cluster-tilting object satisfying
$$C(T') \subset H \subset H_m^{\geq 0}.$$
We obtain that
$$f(T'_i) = f_{T,m}(T'_i) = \langle F({\rm ind}_T(T'_i)), m \rangle \geq 0 $$ for all indecomposable direct summands $T'_i$ of $T$.

Similarly, there exists some cluster-tilting object $T''$ such that $f(T''_i) \leq 0$ for all indecomposable direct
summands $T''_i$ of $T''$.
\end{pf}

\subsection{Another approach to the main theorem}
Let ${\mathcal{C}}_Q$ be the cluster category of a Dynkin quiver
$Q$. In this subsection, we will see another approach to Theorem
\ref{15} by using the work of V. Fock and A. Goncharov \cite{FG}.
For simplicity, we write ${\mathbb{Z}}_{tr}$ for the tropical
semifield $({\mathbb{Z}}, \odot, \oplus)$.

Let ${\mathcal{A}}_{{Q}^{op}} ({\mathbb{Z}}_{tr})$ and
${\mathcal{X}}_{{Q}^{op}} ({\mathbb{Z}}_{tr})$ be the set of
tropical ${\mathbb{Z}}$-points of ${\mathcal{A}}$-variety and
${\mathcal{X}}$-variety \cite{FG} associated with the opposite
quiver $Q^{op}$, respectively. For a vertex $k$ of $Q$, the mutation
$\mu_k: {\mathcal{A}}_{Q^{op}}({\mathbb{Z}}_{tr}) \rightarrow
{\mathcal{A}}_{\mu_k(Q^{op})}({\mathbb{Z}}_{tr})$ is given by the
tropicalization of formula (14) in \cite{FG}:
$$A_k + (\mu_k A)_k = \max \{ \sum_j [b_{jk}]_+ A_j, \sum_j [b_{kj}]_+ A_j \},$$ where $[b_{rs}]_+$ is the number of arrows from $r$
to $s$ in $Q$ (or from $s$ to $r$ in $Q^{op}$). Let $T$ be the image
of $kQ$ in ${\mathcal{C}}_Q$. Then for each tropical
${\mathbb{Z}}$-point $A$ in
${\mathcal{A}}_{Q^{op}}({\mathbb{Z}}_{tr})$, there is a unique
tropical frieze $h$ on ${\mathcal{C}}_Q$ such that $h(T_j) = A_j$
for each $1 \leq j \leq n$. Moreover, this correspondence commutes
with mutation. Besides, we know from \cite{Pla11} that the
isomorphism ${\mathcal{X}}_{Q^{op}}({\mathbb{Z}}_{tr}) \simeq
K_0({\rm add}T)$ commutes with mutation. Given a seed
${\underline{i}}$, in \cite{FG} V. Fock and A. Goncharov considered
the function $P_{\underline{i}} = \sum^n_{i=1} a_i x_i$ on
${\mathcal{A}}({\mathbb{Z}}_{tr}) \times
{\mathcal{X}}({\mathbb{Z}}_{tr})$. Now we can transform the function
$P_{{\underline{i}}}$ in our case as $$P_S = \sum^n_{i=1} h(S_i)
[{\rm ind}_S(Y) : S_i]$$ where $S$ is the cluster-tilting object of
${\mathcal{C}}_Q$ corresponding to the seed ${\underline{i}}$, the
elements $a_i$ correspond to $h(S_i)$ and $x_i$ correspond to $[{\rm
ind}_S(Y) : S_i]$ for some object $Y$ of ${\mathcal{C}}_Q$.

Let $f$ be a tropical frieze on ${\mathcal{C}}_Q$. Let $L$ and $R$
be the same objects as in the proof of Theorem \ref{15}. Assume $X$
is an object of ${\mathcal{C}}_Q$ with
$${\rm ind}_T (X) = [R] - [L]\, (= g(T)).$$ For example, the cone $C(\eta)$ as in the proof of Theorem \ref{15}.
For the pair $N = (f, {\rm ind}_T(X)) \in
{\mathcal{A}}({\mathbb{Z}}_{tr}) \times
{\mathcal{X}}({\mathbb{Z}}_{tr})$, by Theorem 5.2 in \cite{FG},
there exists a cluster-tilting object $T'$ such that all coordinates
$[{\rm ind}_{T'}(X) : T'_i]$ are non-negative. It follows that there
exists some rigid object $X_0 \in {\rm add}T'$ with the same index
as $X$. Set $T'' = {\Sigma}^{-1} T'$, as in the proof of Theorem
\ref{15}, we can also obtain that
$$g(T'') = {\rm ind}_{T''} (X) = {\rm ind}_{T''} (X_0) \in {\mathbb{Z}}^n_{\leq 0}. $$ This gives another approach to the main theorem.

Moreover, our definition for positive cones in Step 1 in the proof of Theorem \ref{17} coincides with Fock-Goncharov's.
From the equality $$K_0 ({\rm add}T) = \bigcup_S C(S),$$ where $S$ ranges over all (finitely many) cluster-tilting objects of ${\mathcal{C}}_Q$, we can also obtain that a finite type cluster ${\mathcal{X}}$-variety
is of definite type (see Corollary 5.5 and Conjecture 5.7 in \cite{FG}).

\vspace{.3cm}

\section{proof of a conjecture of ringel}

\begin{conj} [Section 6, \cite{Rin11}] \label{5}
Let $\Gamma = {\mathbb{Z}} \Delta$ where $\Delta$ is one of the Dynkin diagrams ${\mathbb{A}}_n, {\mathbb{D}}_n, {\mathbb{E}}_6, {\mathbb{E}}_7,
{\mathbb{E}}_8$ and let $f$ be cluster-additive on $\Gamma$. Then $f$ is a non-negative linear combination of cluster-hammock functions (and
therefore of the form
$${\sum}_{x \in {\mathcal{T}}} n_x h_x$$ for a tilting set $\mathcal{T}$ and integers $n_x \in {\mathbb{N}}_0$, for all $x \in {\mathcal{T}}$).
\end{conj}

$\mathbf{Proof \, of \, Conjecture\, \ref{5}}.$
Let $Q$ be an orientation of the Dynkin diagram $\Delta$. Then $\Gamma$ can be viewed as the Auslander-Reiten quiver of the bounded derived category
${\mathcal{D}}_Q$ of the category ${\rm mod}kQ$. Let $I_i$ be the $i$-th indecomposable injective right $kQ$-module. Define a dimension vector
$\underline{d} = (d_i)_{i \in Q_0}$
\begin{displaymath}
d_i = \left\{ \begin{array}{ll}
       f(I_i) & \textrm{if} \, f(I_i) > 0,\\
       0 & \textrm{otherwise}.
\end{array} \right.
\end{displaymath}
Let rep($Q,\underline{d}$) be the affine variety of representations of the opposite quiver $Q^{op}$ with dimension vector $\underline{d}$.
Choose a right $kQ$-module $M$ whose associated point in rep($Q,\underline{d}$) is generic, so that $M$ is rigid.

Define an object $T$ of the cluster category ${\mathcal{C}}_Q$ as $M\oplus( \bigoplus_{f(I_i)<0}(\Sigma P_i)^{-f(I_i)})$. For each $i$ satisfying
$f(I_i) < 0$, we have the following isomorphisms $${\rm Ext}^1_{{\mathcal{C}}_Q} (\Sigma P_i, M) \simeq {\rm Hom}_{{\mathcal{C}}_Q}(P_i, M)
\simeq {\rm Hom}_{kQ}(P_i,M),$$ where the second isomorphism follows from Proposition 1.7 (d) in \cite{BMRRT06}. Notice that the space ${\rm Hom}_{kQ}
(P_i,M)$ vanishes since $M$ does not contain $S_i$ as a composition factor. Thus, the object $T$ is rigid.

Let $M = M_1^{a_1} \oplus \ldots \oplus M_r^{a_r}$ be a decomposition of $M$ with $M_j \,(1 \leq j \leq r)$ indecomposable and pairwise non-isomorphic.
Let $\mathcal{T}$ be the set
$$\{M_j | 1 \leq j \leq r\} \cup \{\Sigma P_i | i \in Q_0 \,\mbox{such that}\, f(I_i)<0\}.$$ Then $\mathcal{T}$ is a partial tilting set \cite{Rin11}. Denote by
$\Sigma {\mathcal{T}}$ the set $\{\Sigma Y | Y \in {\mathcal{T}}\} ( = \{\Sigma M_j | 1 \leq j \leq r \} \cup \{ I_i |  i \in Q_0 \,\mbox{such that}\, f(I_i)<0\})$.
Let $T^{+}$ be a basic cluster-tilting object of ${\mathcal{C}}_Q$ which contains every element in $\mathcal{T}$ as a direct summand. For an indecomposable
object $X$, we use the notation $[N : X]$ to denote the multiplicity of $X$ appearing as a direct summand in ${\mathcal{C}}_Q$ of an object $N$.

Define a new function $f'$ as $\sum_{X \in \Sigma{\mathcal{T}}} [\Sigma T : X] h_X$. Then $f'$ is a cluster-additive function by the Corollary
in Section 6 of \cite{Rin11}. Notice that $$[\Sigma T : I_i] = [T : \Sigma P_i] = -f(I_i)  \quad \mbox{and} \quad
[\Sigma T : \Sigma M_j] = a_j \, (1 \leq j \leq r).$$ Now we rewrite $f'$ as
$$%f' = \sum_{X \in \Sigma{\mathcal{T}}} [\Sigma T : X] h_X =
\sum_{I_i \in \Sigma{\mathcal{T}}} [\Sigma T: I_i] h_{I_i} +
\sum_{\Sigma M_j \in \Sigma{\mathcal{T}}} [\Sigma T : \Sigma M_j] h_{\Sigma M_j} = \sum_{I_i \in \Sigma{\mathcal{T}}} (-f(I_i)) h_{I_i} + \sum^r_{j=1}
a_j h_{\Sigma M_j}.$$

In the following we will show that $f$ and $f'$ coincide on all indecomposable injective $kQ$-modules. Recall that for any pair $X \neq X'$ in a partial
tilting set, the value $h_X(X')$ is zero (Section 5, \cite{Rin11}).

Step 1. Look at the indecomposable injective $kQ$-module $I_l$ satisfying $f(I_l) < 0$.

It is easy to see that
$$f'(I_l) = -f(I_l)h_{I_l} (I_l) = f(I_l).$$

Step 2. Look at the indecomposable injective $kQ$-module $I_l$ satisfying $f(I_l) = 0$.

We have the following isomorphisms
$${\rm Ext}^1_{{\mathcal{C}}_Q}(T, {\Sigma}^{-1}I_l) \simeq {\rm Hom}_{{\mathcal{C}}_Q}(T,I_l) \simeq {\rm Hom}_{{\mathcal{C}}_Q}(T,{\Sigma}^2 P_l)
\simeq D {\rm Hom}_{{\mathcal{C}}_Q}(P_l,T)$$ $$\simeq D {\rm Hom}_{kQ}(P_l,M)
\oplus D {\rm Ext}^1_{{\mathcal{C}}_Q}(P_l,\bigoplus_{f(I_i)<0} (-f(I_i))P_i) = 0.$$
Hence, the set $\Sigma {\mathcal{T}} \cup \{I_l | f(I_l) = 0\}$ is a partial tilting set, which implies that
$$h_X(I_l) = 0, \quad X \in \Sigma {\mathcal{T}}.$$ As a result, we obtain that
$$f'(I_l) = 0 = f(I_l).$$

Step 3. Look at the indecomposable injective $kQ$-module $I_l$ satisfying $f(I_l) > 0$.

We compute the dimension of
${\rm Hom}_{{\mathcal{C}}_Q}(T, I_l)$. As in step 2, we obtain the following isomorphisms
$${\rm Hom}_{{\mathcal{C}}_Q}(T,I_l) \simeq D {\rm Hom}_{kQ}(P_l,M) \simeq {\rm Hom}_{kQ}(M,I_l).$$ It follows that
$${\rm dim Hom}_{{\mathcal{C}}_Q}(T,I_l) = {\rm dim Hom}_{{\mathcal{C}}_Q}(M,I_l) = d_l = f(I_l).$$
Let $B$ denote the endomorphism algebra ${\rm End}_{{\mathcal{C}}_Q}(T^+)$ and $S_{M_j}$ the simple $B$-module which corresponds to
the indecomposable projective $B$-module ${\rm Hom}_{{\mathcal{C}}_Q}(T^+, M_j)$.
For each object $M_j$, we have that
$${\rm dim Hom}_{{\mathcal{C}}_Q}(M_j,I_l) = {\rm dim Hom}_ B ({\rm Hom}_{{\mathcal{C}}_Q}(T^+,M_j), {\rm Hom}_{{\mathcal{C}}_Q}(T^+,I_l))$$
\begin{center}
 = the multiplicity of $S_{M_j}$ as a composition factor of ${\rm Hom}_{{\mathcal{C}}_Q}(T^+,I_l)$ \\
= $h_{\Sigma M_j} (I_l), \quad \quad \quad \quad \quad \quad \quad \quad \quad \quad \quad \quad \quad \quad \quad
\quad \quad \quad \quad \quad \quad \quad \quad \quad \,$
\end{center}
where the last equality appears in the end of the proof of the Lemma in Section 10 of \cite{Rin11}. Since $h_{I_i} (I_l) = 0$ for all $I_i \in
{\Sigma}{\mathcal{T}}$, the following equalities
\begin{align*}
f(I_l)& = {\rm dim Hom}_{{\mathcal{C}}_Q}(M,I_l) = \sum^r_{j=1} a_j {\rm dim Hom}_{{\mathcal{C}}_Q}(M_j,I_l) \\
& =  \sum^r_{j=1} a_j h_{\Sigma M_j} (I_l) = f'(I_l)
\end{align*}
hold.

Therefore, the cluster-additive functions $f$ and $f'$ coincide on all indecomposable injective $kQ$-modules, which implies that $f$ is equal to $f'$.
This completes the proof.

\end{defn}

\end{document}